# A PROBABILISTIC APPROACH TO THE GEOMETRY OF THE $\ell_P^N$-BALL

By Franck Barthe, Olivier Guédon, Shahar Mendelson and Assaf Naor

*Université Paul Sabatier, Université Paris 6, Australian National University and Microsoft Research*

This article investigates, by probabilistic methods, various geometric questions on $B_p^n$, the unit ball of $\ell_p^n$. We propose realizations in terms of independent random variables of several distributions on $B_p^n$, including the normalized volume measure. These representations allow us to unify and extend the known results of the sub-independence of coordinate slabs in $B_p^n$. As another application, we compute moments of linear functionals on $B_p^n$, which gives sharp constants in Khinchine's inequalities on $B_p^n$ and determines the $\psi_2$-constant of all directions on $B_p^n$. We also study the extremal values of several Gaussian averages on sections of $B_p^n$ (including mean width and $\ell$-norm), and derive several monotonicity results as $p$ varies. Applications to balancing vectors in $\ell_2$ and to covering numbers of polyhedra complete the exposition.

**1. Introduction.** For $p > 0$ and a sequence of real numbers $x = (x_i)_{i=1}^{\infty}$ denote $\|x\|_p = (\sum_{i=1}^{\infty} |x_i|^p)^{1/p}$. For $p = \infty$ we set $\|x\|_{\infty} = \sup_{i \in \mathbb{N}} |x_i|$. The space of all infinite sequences $x$ with $\|x\|_p < \infty$ is denoted $\ell_p$. Similarly, the space $\mathbb{R}^n$ equipped with the quasi-norm $\|\cdot\|_p$ is denoted $\ell_p^n$. Finally, the unit balls of $\ell_p^n$ and $\ell_p$ are defined as $B_p^n = \{x \in \mathbb{R}^n; \|x\|_p \leq 1\}$ and $B_p = \{x \in \mathbb{R}^{\mathbb{N}}; \|x\|_p \leq 1\}$, respectively.

The geometry of $\ell_p^n$ spaces in general, and the geometry of the $\ell_p^n$-balls in particular, has been intensively investigated in the past decades. A particular topic of interest has been the evaluation of the extremal volumes of sections and projections of $B_p^n$. Apart from their intrinsic interest, such questions have applications in several probabilistic and geometric contexts, some of which will be described below. The purpose of the present article is to obtain several new results of this flavor. We represent various geometric parameters of $B_p^n$ probabilistically, and apply methods from probability theory to estimate them.

In Section 2 we introduce representations in terms of independent random variables of some distributions on $B_p^n$, including the volume measure on $B_p^n$. Obtaining concrete realizations of the (normalized) volume measure on a general convex body $K \subset \mathbb{R}^n$ seems to be a hopeless task. For general bodies one is therefore reduced to hunting for approximations, and this has been successfully achieved via Markov chain methods by Kannan, Lovasz and Simonovits [18]. (That paper is actually the last in a long list of articles obtaining similar approximate representations. We refer to [18] and the references therein for an accurate historic depiction of the subject.) The simpler structure of $B_p^n$ allows us to give the following representation of the volume measure, which extends to $p > 0$ classical results for $p \in \{1, 2\}$ (see, e.g., Chapter 2 in [15]).

THEOREM 1. *Let $g_1, \ldots, g_n$ be i.i.d. random variables with density $1/(2\Gamma(1+1/p))e^{-|t|^p}$ ($t \in \mathbb{R}$), and let $Z$ be an exponential random variable independent of $g_1, \ldots, g_n$ (i.e., the density of $Z$ is $e^{-t}$, $t \geq 0$). Denote $G = (g_1, \ldots, g_n) \in \mathbb{R}^n$ and consider the random vector*

$$V = \frac{G}{(\sum_{i=1}^n |g_i|^p + Z)^{1/p}}.$$

*Then $V$ generates the normalized volume measure on $B_p^n$, that is, for every measurable $A \subset \mathbb{R}^n$,*

$$P(V \in A) = \frac{\mathrm{vol}(A \cap B_p^n)}{\mathrm{vol}(B_p^n)}.$$

Section 2.1 provides a simple probabilistic perspective to the sub-independence of coordinate slabs on $B_p^n$. This remarkable fact was originally proved by Ball and Perissinaki [4] for the volume measure and in [24] for the cone measure. We establish this property for more general distributions, combining an extension of Theorem 1 with arguments similar to the proof of the classical FKG inequality [16].

In Section 2.2, Theorem 1 is applied to the study of the moments of linear functionals on $B_p^n$ for $p \geq 1$. Answering a question posed to us by Giannopoulos, we estimate the best constants in the Khinchine inequality on $B_p^n$ and describe the so-called $\psi_2$-directions of $B_p^n$.

Section 3 is devoted to the analysis of the extremal values of several geometric parameters of sections of $B_p^n$ for $p > 0$. A classical result of Meyer and Pajor [22] states that for every $k$-dimensional subspace $E$ of $\mathbb{R}^n$, if $p \leq 2$, then $\mathrm{vol}_k(E \cap B_p^n) \leq \mathrm{vol}_k(B_p^k)$, and if $p \geq 2$, then $\mathrm{vol}_k(E \cap B_p^n) \geq \mathrm{vol}_k(B_p^k)$.

More results on critical sections of $B_p^n$ appear in the papers [2, 3, 20, 22], which rely on harmonic analysis methods. In Section 3.1 we show that for every $0 \leq \alpha \leq k$, every $0 \leq \beta \leq p$ and every $k$-dimensional subspace $E$ of $\mathbb{R}^n$, if $0 < p \leq 2$, then

$$\int_{S^{n-1} \cap E} \|x\|_p^{-\alpha} \, dx \leq \int_{S^{k-1}} \|x\|_p^{-\alpha} \, dx \quad \text{and} \quad \int_{S^{n-1} \cap E} \|x\|_p^{\beta} \, dx \geq \int_{S^{k-1}} \|x\|_p^{\beta} \, dx,$$

and if $2 < p \leq \infty$, then

$$\int_{S^{n-1} \cap E} \|x\|_p^{-\alpha} \, dx \geq \int_{S^{k-1}} \|x\|_p^{-\alpha} \, dx \quad \text{and} \quad \int_{S^{n-1} \cap E} \|x\|_p^{\beta} \, dx \leq \int_{S^{k-1}} \|x\|_p^{\beta} \, dx.$$

The case $\alpha = k$ in the above inequalities is just a restatement of the Meyer–Pajor theorem. The case $\beta = p$ follows from the following stronger monotonicity result, proved in Section 3.1, that the mapping

$$p > 0 \mapsto \frac{\int_{S^{n-1} \cap E} \|x\|_p^p \, dx}{\int_{S^{k-1}} \|x\|_p^p \, dx}$$

is increasing in $p$.

Since Gaussian and spherical averages of homogeneous functions are proportional, these facts can be restated in terms of moments of Gaussian vectors. Note that the above quantities encompass useful classical parameters of the geometry of Banach spaces, such as mean width and $\ell$-norm (see, e.g., [29], page 35).

The proofs appear in Section 3.1 and consist of finding probabilistic expressions of various expectations of Gaussian vectors on subspaces of $\mathbb{R}^n$, and then applying stochastic orderings to estimate them.

In Section 3.2 we apply the Brascamp–Lieb inequality to obtain estimates in the other direction.

Section 3.3 deals with the case of the cube $B_\infty^n$. We derive the following distributional inequalities, valid for all $k$-dimensional subspaces $E \subset \mathbb{R}^n$ and every $r > 0$:

$$\gamma_k(rB_\infty^k) \leq \gamma_E(E \cap rB_\infty^n) \leq \gamma_k\left(r\sqrt{\frac{n}{k}} B_\infty^k\right),$$

where $\gamma_k, \gamma_E$ denote the standard Gaussian measure on $\mathbb{R}^k$ and $E$, respectively. The right-hand side of the above inequality follows from the Brascamp–Lieb inequality, and the left-hand side from the following monotonicity result: for every $k$-dimensional subspace $E \subset \mathbb{R}^n$, the function

$$r > 0 \mapsto \frac{\gamma_E(E \cap rB_\infty^n)}{\gamma_k(rB_\infty^k)}$$

is nonincreasing.

Sections 3.4 and 3.5 are devoted to applications of the previous results. Section 3.4 deals with the Komlós conjecture which asks whether there is a universal constant $c > 0$ such that for every $x_1, \ldots, x_m \in B_2^n$, there are signs $\varepsilon_1, \ldots, \varepsilon_m \in \{-1, 1\}$ for which $\|\sum_{i=1}^m \varepsilon_i x_i\|_\infty \leq c$. This challenging problem remains unsolved, and the best upper bound on $c$, due to Banaszczyk [5], is $c = O(\sqrt{\log n})$. We show that our estimates, together with Banaszczyk's theorem, yield an infinite-dimensional version of this result, which implies in particular a better upper bound when $m = o(n)$.

PROPOSITION 1. *There is an absolute constant $C > 0$ such that for every integer $m > 0$ and every $x_1, \ldots, x_m \in \ell_\infty$, there are signs $\varepsilon_1, \ldots, \varepsilon_m \in \{-1, 1\}$ for which*

$$\left\| \sum_{i=1}^m \varepsilon_i x_i \right\|_\infty \leq C\sqrt{\log d} \cdot \max_{1 \leq i \leq m} \|x_i\|_2 \leq C\sqrt{\log m} \cdot \max_{1 \leq i \leq m} \|x_i\|_2,$$

*where $d$ is the dimension of the linear span of $x_1, \ldots, x_m$.*

Section 3.5 answers a question posed to us by Talagrand, concerning the number of cubes required to cover a convex hull of a finite number of points in $\ell_2$. Given two convex sets $K, L \subset \ell_\infty$, denote by $N(K, L)$ the minimal number of translates of $L$ required to cover $K$ (this number may be infinite). Obtaining sharp bounds on this parameter is of fundamental importance in several problems in convex geometry (see, e.g., [29]), probability (see, e.g., [21]) and operator theory (see, e.g., [28]). Given $A \subset \ell_\infty$, we denote by absconv$(A)$ the convex hull of $A \cup (-A)$. The main result of Section 3.5 is:

PROPOSITION 2. *There exists an absolute constant $C > 0$ such that for every integer $m$, $\varepsilon > 0$ and $2 \leq p \leq \infty$, for all $x_1, \ldots, x_m$ in the unit ball of $\ell_2$,*

$$\log N(\mathrm{absconv}\{x_1, \ldots, x_m\}, \varepsilon B_p) \leq C \frac{\log m}{\varepsilon^{p/(p-1)}}.$$

Such a statement is already known for $p = 2$ by the results of Carl and Pajor [14]. From Schütt's results [34] on the entropy of the identity operator between $\ell_2^d$ and $\ell_p^d$, if the points $x_1, \ldots, x_m$ are assumed to be in an ambient $\ell_\infty^d$, then such an inequality is valid with the term $\log m$ replaced by $\log \max(m, d)$. Proposition 2 bounds the covering number of the polyhedron absconv$\{x_1, \ldots, x_m\} \subset B_2^d$ in terms of the number of its vertices, independently of the ambient dimension.

**2. Representation of measures on $B_p^n$.** We begin by stating a probabilistic representation of the cone measure on $\partial B_p^n$ which is due to Schechtman and Zinn [32] and independently to Rachev and Rüschendorf [31].

This representation has applications of probabilistic and geometric nature [9, 24, 25, 33].

Let $K$ be a convex symmetric body in $\mathbb{R}^n$. Recall that the cone measure on $\partial K$, denoted $\mu_K$, is defined for $A \subset \partial K$ by

$$\mu_K(A) = \frac{\operatorname{vol}(ta; a \in A, 0 \le t \le 1)}{\operatorname{vol}(K)}.$$

Thus, $\mu_K(A)$ is the volume of the cone with base $A$ and cusp $0$, normalized by the volume of $K$. Alternately, $\mu_K$ is the unique measure for which the following polar integration formula holds: for every $f \in L_1(\mathbb{R}^n)$,

$$\int_{\mathbb{R}^n} f(x)\, dx = n \cdot \operatorname{vol}(K) \int_0^\infty r^{n-1} \int_{\partial K} f(rz)\, d\mu_K(z)\, dr.$$

Schechtman and Zinn and Rachev and Rüschendorf proved the following.

THEOREM 2 ([31, 32]). *Let $g_1, \ldots, g_n$ be i.i.d. random variables with density $e^{-|t|^p}/(2\Gamma(1+1/p))$, $t \in \mathbb{R}$. Consider the random vector $G = (g_1, \ldots, g_n) \in \mathbb{R}^n$, and denote*

$$Y = \frac{G}{\|G\|_p} = \frac{G}{(\sum_{i=1}^n |g_i|^p)^{1/p}}.$$

*Then $Y$ is independent of $\|G\|_p$. Moreover, $Y$ generates the measure $\mu_{B_p^n}$; that is, for every measurable $A \subset \partial B_p^n$, $\mu_{B_p^n}(A) = P(Y \in A)$.*

We propose the following extension:

THEOREM 3. *Let $G = (g_1, \ldots, g_n)$ be a random vector as in Theorem 2. Let $W$ be a nonnegative random variable with distribution $h$, and independent of $G$. Then the random vector*

$$\frac{G}{(\|G\|_p^p + W)^{1/p}}$$

*generates the measure $h(\{0\})\mu_{B_p^n} + \Psi \lambda_{B_p^n}$, where $\lambda_{B_p^n}$ stands for Lebesgue's measure restricted to $B_p^n$, and for $x \in B_p^n$, $\Psi(x) = \psi(\|x\|_p)$, where for $r \in [0,1]$*

$$(1) \quad \left[\Gamma\left(1 + \frac{1}{p}\right)\right]^n \psi(r) = \frac{1}{(1-r^p)^{n/p+1}} \int_{(0,\infty)} w^{n/p} e^{-r^p w/(1-r^p)}\, dh(w).$$

PROOF. Note that the density of $|g_i|^p$ is

$$\frac{d}{du} P(|g_i| \le u^{1/p}) = \frac{2u^{1/p-1}}{p} \cdot \frac{1}{2\Gamma(1+1/p)} e^{-u}$$

$$= \frac{1}{\Gamma(1/p)} u^{1/p-1} e^{-u}, \qquad u > 0.$$

In other words, $|g_i|^p$ has a $gamma(1/p, 1)$ distribution. By the additivity property of the gamma semigroup, the random variable $\|G\|_p^p = \sum_{i=1}^n |g_i|^p$ has a $gamma(n/p, 1)$ distribution, that is, its density is $1/\Gamma(n/p) u^{n/p-1} e^{-u}$ ($u \geq 0$).

For any $f \in L_1(\mathbb{R}^n)$, and conditioning on $W$,

$$\mathbb{E} f\left(\frac{G}{(\|G\|_p^p + W)^{1/p}}\right) = \int_{[0,\infty)} \mathbb{E} f\left(\frac{G}{(\|G\|_p^p + w)^{1/p}}\right) dh(w).$$

Since $G/\|G\|_p$ and $\|G\|_p$ are independent, then for every $w > 0$,

$$\mathbb{E} f\left(\frac{G}{(\|G\|_p^p + w)^{1/p}}\right)$$
$$= \mathbb{E} f\left(\left(\frac{\|G\|_p^p}{\|G\|_p^p + w}\right)^{1/p} \frac{G}{\|G\|_p}\right)$$
$$= \frac{1}{\Gamma(n/p)} \int_0^\infty u^{n/p-1} e^{-u} \mathbb{E} f\left(\left(\frac{u}{u+w}\right)^{1/p} \frac{G}{\|G\|_p}\right) du$$
$$= \frac{1}{\Gamma(n/p)} \int_0^1 \left(\frac{r^p w}{1-r^p}\right)^{n/p-1} e^{-r^p w/(1-r^p)} \cdot \mathbb{E} f\left(r \frac{G}{\|G\|_p}\right) \frac{p r^{p-1} w}{(1-r^p)^2} dr,$$

where we have made the change of variable $\frac{u}{u+w} = r^p$. Hence,

$$\mathbb{E} f\left(\frac{G}{(\|G\|_p^p + W)^{1/p}}\right) - h(\{0\}) \mathbb{E} f\left(\frac{G}{\|G\|_p}\right)$$
$$= \frac{p}{\Gamma(n/p)} \int_{(0,\infty)} w^{n/p}$$
$$\qquad \times \int_0^1 \frac{r^{n-1}}{(1-r^p)^{n/p+1}} e^{-r^p w/(1-r^p)} \cdot \mathbb{E} f\left(r \frac{G}{\|G\|_p}\right) dr \, dh(w)$$
$$= \frac{n}{\Gamma(n/p+1)} \int_0^1 \frac{r^{n-1}}{(1-r^p)^{n/p+1}}$$
$$\qquad \times \left(\int_{(0,\infty)} w^{n/p} e^{-r^p w/(1-r^p)} dh(w)\right) \mathbb{E} f\left(r \frac{G}{\|G\|_p}\right) dr.$$

On the other hand, let $M$ be a probability measure on $B_p^n$ with $\ell_p$-radial density $\phi(\|x\|_p)$ ($x \in B_p^n$). By the polar coordinate integration formula for $\mu_{B_p^n}$, the representation from Theorem 2 and the fact that $\mathrm{vol}(B_p^n) = \frac{[2\Gamma(1/p+1)]^n}{\Gamma(n/p+1)}$ (see, e.g., page 11 in [29]),

$$\int_{\mathbb{R}^n} f(x) \, dM(x) = n \, \mathrm{vol}(B_p^n) \int_0^1 r^{n-1} \phi(r) \mathbb{E} f\left(r \frac{G}{\|G\|_p}\right) dr$$

$$= \frac{n[2\Gamma(1/p+1)]^n}{\Gamma(n/p+1)} \int_0^1 r^{n-1} \phi(r) \mathbb{E} f\left(r \frac{G}{\|G\|_p}\right) dr,$$

from which the result easily follows. □

Since (1) holds true for

$$dh(w) = e^{-w} \mathbf{1}_{\{w>0\}} \, dw \quad \text{and} \quad \psi(r) = \frac{\mathbf{1}_{[0,1]}(r)}{\text{vol}(B_p^n)},$$

we have established Theorem 1. We now study more general distributions. By making the change of variable $s = \frac{r^p}{1-r^p}$ in (1), we obtain the following representation theorem. We refer to [36] for completely monotone functions and the Laplace transform.

THEOREM 4. *Let $\nu$ be a probability measure on $\mathbb{R}^n$ with density $\psi(\|x\|_p)\mathbf{1}_{[0,1]}(\|x\|_p)$. Assume that the function*

$$s \mapsto \frac{1}{(1+s)^{n/p+1}} \psi\left(\left(\frac{s}{1+s}\right)^{1/p}\right), \qquad s > 0,$$

*is completely monotone. Then there is a positive random variable $W$ such that for every measurable $A \subset \mathbb{R}^n$,*

$$\nu(A) = P\left(\frac{G}{(\|G\|_p^p + W)^{1/p}} \in A\right),$$

*and the density of $W$ is given by*

$$\frac{[2\Gamma(1/p+1)]^n}{w^{n/p}} \mathcal{L}^{-1}\left[s \mapsto \frac{1}{(1+s)^{n/p+1}} \psi\left(\left(\frac{s}{1+s}\right)^{1/p}\right)\right](w), \qquad w > 0,$$

*where $\mathcal{L}$ is the Laplace transform.*

Next, we single out an interesting case for which the above theorem may be applied: when $W$ is a *gamma$(\alpha, 1)$* random variable, the density of $W$ is $h(w) = 1/\Gamma(\alpha) w^{\alpha-1} e^{-w}$, and thus

$$\left[2\Gamma\left(\frac{1}{p}+1\right)\right]^n \psi(r) = \frac{1}{\Gamma(\alpha)(1-r^p)^{n/p+1}} \int_0^\infty w^{n/p+\alpha-1} e^{-w/(1-r^p)} \, dw$$

$$= \frac{(1-r^p)^{n/p+\alpha}}{\Gamma(\alpha)(1-r^p)^{n/p+1}} \int_0^\infty w^{n/p+\alpha-1} e^{-w} \, dw$$

$$= \frac{(1-r^p)^{\alpha-1} \Gamma(n/p+\alpha)}{\Gamma(\alpha)}.$$

COROLLARY 3. *Let $W$ be a gamma$(\alpha, 1)$ random variable. Then the random vector $\frac{G}{(\|G\|_p^p + W)^{1/p}}$ generates the measure on $B_p^n$ with density*

(2) $$f(x) = \frac{\Gamma(n/p + \alpha)}{\Gamma(\alpha)[2\Gamma(1/p+1)]^n}(1 - \|x\|_p^p)^{\alpha-1}\mathbf{1}_{[0,1]}(\|x\|_p).$$

Finally let us give a geometric interpretation of some of our representations. Fix two integers $m, n$ and consider the orthogonal projection of the cone measure on $\partial B_p^{n+m}$ onto the first $n$ coordinates. By the Schechtman–Zinn theorem, this measure is generated by the random vector

$$\frac{(g_1, \ldots, g_n)}{(\sum_{i=1}^n |g_i|^p + \sum_{i=n+1}^{m+n} |g_i|^p)^{1/p}}.$$

The random variable $\sum_{i=n+1}^{m+n} |g_i|^p$ is independent of $g_1, \ldots, g_n$ and has a gamma$(m/p, 1)$ distribution. Hence, the above discussion leads to the following extension of classical observations about $B_1^n$ and $B_2^n$ (for these sets the cone measure coincides with the better studied normalized surface measure).

COROLLARY 4. *When $p$ is an integer, the orthogonal projection of the cone measure on $\partial B_p^{n+p}$ onto the first $n$ coordinates is the (normalized) volume measure on $B_p^n$. More generally, for arbitrary $p > 0$, the orthogonal projection of the cone measure on $\partial B_p^{n+m}$ onto the first $n$ coordinates has density*

$$f(x) = \frac{\Gamma((n+m)/p)}{\Gamma(m/p)[2\Gamma(1/p+1)]^n}(1 - \|x\|_p^p)^{m/p-1}\mathbf{1}_{[0,1]}(\|x\|_p).$$

2.1. *An application*: *sub-independence of coordinate slabs.* The sub-independence of coordinate slabs in $B_p^n$ is helpful in the study of the central limit problem [1, 25] and of various deviation inequalities [8, 24]. More precisely, this property is enjoyed by the normalized volume measure on $B_p^n$, as proved analytically in [4] and geometrically in [1]. It was established probabilistically in [24] for the cone measure on $B_p^n$. In this section we combine our representation results with an argument of [24] in order to derive sub-independence of coordinate slabs for a wider class of distributions. We require the following result:

THEOREM 5 ([7]). *Let $X_1, \ldots, X_n$ be independent symmetric random variables. Assume that $X_i$ has density $\psi_i = e^{-V_i}$, where $V_i$ is locally integrable. For $X = (X_1, \ldots, X_n)$, the random vector $\frac{X}{\|X\|_p}$ is independent of the random variable $\|X\|_p$ if and only if there are $b_1, \ldots, b_n > -1$ and $a, c_1, \ldots, c_n > 0$ such that for every $1 \leq i \leq n$, $\psi_i(x) = c_i|x|^{b_i}e^{-a|x|^p}$.*

REMARK. As a consequence of this characterization, setting for $k \leq n$, $X^k := (X_1, \ldots, X_k)$ (where we write for simplicity $X$ for $X^n$), it follows that the independence of $\frac{X}{\|X\|_p}$ from $\|X\|_p$ guarantees for every $k < n$ the independence of $\frac{X^k}{\|X^k\|_p}$ from $\|X^k\|_p$.

The following lemma was essentially proved in [24]. It was stated there for the cone measure on $\partial B_p^n$, but the proof carries through to the more general setting. We sketch the argument for the sake of completeness. Our geometric interest led us to consider symmetric variables, but it is clear that the result concerns nonnegative variables.

LEMMA 5. *Let $X_1, \ldots, X_n$ be independent symmetric random variables. For $i = 1, \ldots, n-1$, assume that $X_i$ has density $\psi_i = \exp(-V_i)$, where $V_i$ is locally integrable. We write $\mu_n$ for the law of $|X_n|$. Denote $X = (X_1, \ldots, X_n)$, $X^{n-1} = (X_1, \ldots, X_{n-1})$ and assume that $\frac{X^{n-1}}{\|X^{n-1}\|_p}$ is independent of $\|X^{n-1}\|_p$. Let $f_1, \ldots, f_n : [0, \infty) \to [0, \infty)$ be nonnegative nondecreasing functions. Then*

$$\mathbb{E}\left[\prod_{i=1}^n f_i\left(\frac{|X_i|}{\|X\|_p}\right)\right] \leq \prod_{i=1}^n \mathbb{E} f_i\left(\frac{|X_i|}{\|X\|_p}\right).$$

PROOF. The proof is by induction on $n$. Assume that $n > 1$ and that the required inequality holds for $n - 1$. Conditioning on $|X_n|$,

$$\mathbb{E}\left[\prod_{i=1}^n f_i\left(\frac{|X_i|}{\|X\|_p}\right)\right] = \int_{\mathbb{R}^+} \mathbb{E}\left\{\left[\prod_{i=1}^{n-1} f_i\left(\frac{|X_i|}{(\|X^{n-1}\|_p^p + r^p)^{1/p}}\right)\right] \right.$$
$$\left. \times f_n\left(\frac{r}{(\|X^{n-1}\|_p^p + r^p)^{1/p}}\right)\right\} d\mu_n(r).$$

Note that by the remark after Theorem 5, $\frac{X^{n-2}}{\|X^{n-2}\|_p}$ and $\|X^{n-2}\|_p$ are independent, so that we may apply the inductive hypothesis. Denote by $\varphi$ the density of $\|X^{n-1}\|_p$, and by the independence of $\frac{X^{n-1}}{\|X^{n-1}\|_p}$, and $\|X^{n-1}\|_p$ it follows that for every $r > 0$,

$$\mathbb{E}\left\{\left[\prod_{i=1}^{n-1} f_i\left(\frac{|X_i|}{(\|X^{n-1}\|_p^p + r^p)^{1/p}}\right)\right] \cdot f_n\left(\frac{r}{(\|X^{n-1}\|_p^p + r^p)^{1/p}}\right)\right\}$$

$$= \int_0^\infty \varphi(u) f_n\left(\frac{r}{(u^p + r^p)^{1/p}}\right) \cdot \mathbb{E}\left[\prod_{i=1}^{n-1} f_i\left(\frac{u}{(u^p + r^p)^{1/p}} \cdot \frac{|X_i|}{\|X^{n-1}\|_p}\right)\right] du$$

$$\leq \int_0^\infty \varphi(u) f_n\left(\frac{r}{(u^p + r^p)^{1/p}}\right) \cdot \prod_{i=1}^{n-1} \mathbb{E} f_i\left(\frac{u}{(u^p + r^p)^{1/p}} \cdot \frac{|X_i|}{\|X^{n-1}\|_p}\right) du.$$

For $u>0$ let $h_u(r) = f_n(\frac{r}{(u^p+r^p)^{1/p}})$ and

$$k_u(r) = \prod_{i=1}^{n-1} \mathbb{E} f_i\left(\frac{u}{(u^p+r^p)^{1/p}} \cdot \frac{|X_i|}{\|X^{n-1}\|_p}\right).$$

Thus $h_u$ is nondecreasing and $k_u$ is nonincreasing and if $X'_n$ is an independent copy of $X_n$, then $[h_u(|X_n|) - h_u(|X'_n|)] \cdot [k_u(|X_n|) - k_u(|X'_n|)] \leq 0$ pointwise. Taking expectation of this inequality,

$$\int_{\mathbb{R}^+} h_u(r) k_u(r) \, d\mu_n(r) \leq \left(\int_{\mathbb{R}^+} h_u(r) \, d\mu_n(r)\right)\left(\int_{\mathbb{R}^+} k_u(r) \, d\mu_n(r)\right),$$

implying that

$$\mathbb{E}\left[\prod_{i=1}^n f_i\left(\frac{|X_i|}{\|X\|_p}\right)\right] \leq \int_{\mathbb{R}^+} \int_0^\infty \varphi(u) h_u(r) k_u(r) \, du \, d\mu_n(r)$$

$$\leq \int_0^\infty \varphi(u)\left(\int_{\mathbb{R}^+} h_u(r) \, d\mu_n(r)\right)\left(\int_{\mathbb{R}^+} k_u(r) \, d\mu_n(r)\right) du$$

$$= \prod_{i=1}^n \mathbb{E} f_i\left(\frac{|X_i|}{\|X\|_p}\right). \qquad \square$$

The main result of this section is contained in the following theorem.

THEOREM 6. *Let $G = (g_1, \ldots, g_n)$ be a random vector with independent coordinates with distribution $e^{-|t|^p}/(2\Gamma(1+1/p))$, $t \in \mathbb{R}$. Let $W$ be a non-negative random variable, independent from $G$. Let $\nu$ be the distribution (supported on $B_p^n$) of the vector*

$$\frac{G}{(\|G\|_p^p + W)^{1/p}}.$$

*Then for every $s_1, \ldots, s_n > 0$,*

$$\nu\left(\bigcap_{i=1}^n \{|x_i| \geq s_i\}\right) \leq \prod_{i=1}^n \nu(\{|x_i| \geq s_i\}).$$

PROOF. Assume that $\varepsilon$ is a random variable independent of $G$ and $W$ which takes the values $+1, -1$ with probability $1/2$. We set $X = (g_1, \ldots, g_n, \varepsilon W^{1/p}) \in \mathbb{R}^{n+1}$. By Theorem 2, $\frac{G}{\|G\|_p}$ and $\|G\|_p$ are independent, so we can apply Lemma 5 to $X$, with $f_i(x) = \mathbf{1}_{[s_i,\infty)}(x)$ for $i = 1, \ldots, n$ and $f_{n+1} = 1$. Hence,

$$P\left(\bigcap_{i=1}^n \left\{\frac{|g_i|}{(\|G\|_p^p + W)^{1/p}} \geq s_i\right\}\right) \leq \prod_{i=1}^n P\left(\left\{\frac{|g_i|}{(\|G\|_p^p + W)^{1/p}} \geq s_i\right\}\right). \qquad \square$$

REMARK. By the very same proof, one can see that the conclusion of Lemma 5 holds for nonnegative, nonincreasing functions. Thus Theorem 6 also holds for symmetric slabs $\{|x_i| \leq s_i\}$.

REMARK. We have obtained sub-independence of coordinate slabs for a class of measures on $B_p^n$, described in Theorem 3. This unifies the previously known occurrences of such sub-independence, since the cone measure $\mu_p^n$ and the normalized volume measure on $B_p^n$ belong to this class. We obtain new concrete examples, as the measures $\nu_\alpha$ with density

$$f_\alpha(x) = \frac{\Gamma(n/p+\alpha)}{\Gamma(\alpha)[2\Gamma(1/p+1)]^n}(1-\|x\|_p^p)^{\alpha-1}\mathbf{1}_{[0,1]}(\|x\|_p).$$

Since these measures $\nu_\alpha$ are isotropic, an immediate consequence of Theorem 6 is that they enjoy the central limit property in the sense that Theorem 5 of [25] holds for them. We refer to that paper for details.

2.2. *An application*: *moment inequalities on $B_p^n$ for $p \geq 1$.* In what follows, given two sequences of positive real numbers $(a_i)_{i \in I}, (b_i)_{i \in I}$, the notation $a_i \sim b_i$ refers to the fact that there are constants $c$ and $C$ such that for all $i \in I$, $ca_i \leq b_i \leq Ca_i$. We emphasize that such $c, C$ are always absolute numerical constants.

We can relate moments of linear functionals on $B_p^n$ to moments of linear functionals of the random vector $G = (g_1, \ldots, g_n)$ with independent coordinates with distribution $e^{-|t|^p}/(2\Gamma(1+1/p))$:

LEMMA 6. *For every integer $n \geq 1$, every $p, q \geq 1$ and every $a \in \mathbb{R}^n$, one has*

$$\left(\frac{1}{\mathrm{vol}(B_p^n)}\int_{B_p^n}\left|\sum_{i=1}^n a_i x_i\right|^q dx\right)^{1/q} \sim \frac{1}{(\max\{n,q\})^{1/p}}\left(\mathbb{E}\left|\sum_{i=1}^n a_i g_i\right|^q\right)^{1/q}.$$

PROOF. Denote $a = (a_1, \ldots, a_n)$. By the probabilistic representation of the volume measure on $B_p^n$ established in Theorem 1,

$$\frac{1}{\mathrm{vol}(B_p^n)}\int_{B_p^n}\left|\sum_{i=1}^n a_i x_i\right|^q dx = \mathbb{E}\left|\left\langle\frac{G}{(\|G\|_p^p+Z)^{1/p}},a\right\rangle\right|^q$$

$$= \mathbb{E}\left[\left(\frac{\|G\|_p^p}{\|G\|_p^p+Z}\right)^{q/p}\left|\left\langle\frac{G}{\|G\|_p},a\right\rangle\right|^q\right]$$

$$= \left[\mathbb{E}\left(\frac{\|G\|_p^p}{\|G\|_p^p+Z}\right)^{q/p}\right]\cdot\left[\mathbb{E}\left|\left\langle\frac{G}{\|G\|_p},a\right\rangle\right|^q\right]$$

$$= \left[\mathbb{E}\left(\frac{\|G\|_p^p}{\|G\|_p^p+Z}\right)^{q/p}\right]\cdot\frac{\mathbb{E}|\langle G,a\rangle|^q}{\mathbb{E}\|G\|_p^q},$$

where we have used the independence of $\frac{G}{\|G\|_p}$ and $\|G\|_p$. Applying this identity to $a = (1, 0, \ldots, 0)$ yields

$$\frac{1}{\mathbb{E}\|G\|_p^q}\left[\mathbb{E}\left(\frac{\|G\|_p^p}{\|G\|_p^p + Z}\right)^{q/p}\right] = \frac{1}{\text{vol}(B_p^n)\mathbb{E}|g_1|^q}\int_{B_p^n}|x_1|^q\,dx.$$

Now, $\mathbb{E}|g_1|^q = \frac{\Gamma((q+1)/p+1)}{(q+1)\Gamma(1/p+1)}$, and for every $p, q \geq 1$,

$$\frac{1}{\text{vol}(B_p^n)}\int_{B_p^n}|x_1|^q\,dx$$

$$= \frac{2\,\text{vol}(B_p^{n-1})}{\text{vol}(B_p^n)}\int_0^1 u^q(1-u^p)^{(n-1)/p}\,du$$

$$= \frac{2[2\Gamma(1/p+1)]^{n-1}\Gamma(n/p+1)}{\Gamma((n-1)/p+1)[2\Gamma(1/p+1)]^n}\frac{1}{p}\int_0^1 v^{(q+1)/p-1}(1-v)^{(n-1)/p}\,dv$$

$$= \frac{\Gamma(n/p+1)}{\Gamma((n-1)/p+1)\Gamma(1/p+1)}\cdot\frac{\Gamma((q+1)/p+1)\Gamma((n-1)/p+1)}{(q+1)\Gamma((n+q)/p+1)},$$

where we have used $\text{vol}(B_p^n) = (\Gamma(1+1/p))^n/\Gamma(1+n/p)$. Therefore,

$$\frac{1}{\mathbb{E}\|G\|_p^q}\left[\mathbb{E}\left(\frac{\|G\|_p^p}{\|G\|_p^p + Z}\right)^{q/p}\right] = \frac{\Gamma(n/p+1)}{\Gamma((n+q)/p+1)},$$

and by Stirling's formula, there are constants $c, C > 0$ such that for all $n, q, p \geq 1$,

$$c\frac{1}{(\max\{n,q\})^{1/p}} \leq \left(\frac{\Gamma(n/p+1)}{\Gamma((n+q)/p+1)}\right)^{1/q} \leq C\frac{1}{(\max\{n,q\})^{1/p}}. \qquad \square$$

For independent symmetric random variable with log-concave cumulated distribution function, Gluskin and Kwapień [17] obtained an almost exact expression of moments of linear functionals. We apply their result to obtain:

PROPOSITION 7.  *Let $n \geq 1$ be an integer. Let $p, q \geq 1$ and $a_1 \geq a_2 \geq \cdots \geq a_n \geq 0$. Then*

$$\left(\mathbb{E}\left|\sum_{i=1}^n a_i g_i\right|^q\right)^{1/q} \sim q^{1/p}\|(a_i)_{i\leq q}\|_{p'} + \sqrt{q}\|(a_i)_{i>q}\|_2,$$

*where $p' \in [1, +\infty]$ is the dual exponent of $p$, defined by $\frac{1}{p} + \frac{1}{p'} = 1$.*

The proof of Proposition 7 requires some preparation.

LEMMA 8. *For every $t > 0$,*
$$\int_t^\infty e^{-u^p} du \leq \frac{e^{-t^p}}{pt^{p-1}},$$

*and for every $t \geq 1$,*
$$\int_t^\infty e^{-u^p} du \geq \frac{e^{-t^p}}{2pt^{p-1}}.$$

*In addition, the function $t \mapsto \int_t^\infty e^{-u^p} du$ is log-concave.*

PROOF. For every $t > 0$,
$$\int_t^\infty e^{-u^p} du \leq \int_t^\infty \frac{u^{p-1}}{t^{p-1}} e^{-u^p} du = \frac{e^{-t^p}}{pt^{p-1}}.$$

To prove the reverse inequality assume that $t \geq 1$. Integrating by parts,
$$\int_t^\infty e^{-u^p} du = \int_t^\infty u^{1-p} \cdot u^{p-1} e^{-u^p} du$$
$$= \frac{e^{-t^p}}{pt^{p-1}} - \frac{p-1}{p} \int_t^\infty \frac{e^{-u^p}}{u^p} du \geq \frac{e^{-t^p}}{pt^{p-1}} - \int_t^\infty e^{-u^p} du,$$

which implies the assertion.

Finally, set $f(t) = \int_t^\infty e^{-u^p} du$. In order to show that $f$ is log-concave it suffices to show that $f''f - (f')^2 \leq 0$ point-wise. Now,
$$f''(t)f(t) - f'(t)^2 = e^{-t^p}\left(pt^{p-1} \int_t^\infty e^{-u^p} du - e^{-t^p}\right) \leq 0,$$

by the first assertion we proved. □

PROOF OF PROPOSITION 7. In what follows $g$ denotes a random variable with density $1/(2\Gamma(1+1/p))e^{-|t|^p}$. Let $\theta_p > 0$ be such that $P(\theta_p |g| \geq 1) = 1/e$. Denote $N(t) = -\log P(\theta_p |g| \geq t)$ and let $N^*(t)$ be the Legendre transform of $N$, that is, $N^*(t) = \sup\{ts - N(s); s > 0\}$. By Lemma 8, $N$ is convex, and a result of Gluskin and Kwapień [17] states that in this case,
$$\left(\mathbb{E}\left|\sum_{i=1}^n a_i g_i\right|^q\right)^{1/q} \sim \theta_p \left[\inf\left\{t > 0; \sum_{i \leq q} N^*\left(\frac{qa_i}{t}\right) \leq q\right\} + \sqrt{q}\left(\sum_{i > q} a_i^2\right)^{1/2}\right],$$

where $a_1 \geq a_2 \geq \cdots \geq a_n \geq 0$. When $p = 1$, all the above quantities are easily computed [in particular $N(t) = t$] and the proposition follows. For $p > 1$, we shall prove below that there exist universal constants $c, c', C, C' > 0$ such that for all $p > 1$,

(3)
$$c' \leq \theta_p \leq C', \quad \forall t > 0 \quad (N^*(t))^{(p-1)/p} \leq Ct \quad \text{and}$$
$$\forall t \geq 2 \quad (N^*(t))^{(p-1)/p} \geq ct.$$

First we explain how these inequalities allow us to conclude. Let

$$t_0 = \inf\left\{t > 0; \sum_{i \leq q} N^*\left(\frac{qa_i}{t}\right) \leq q\right\}.$$

The above upper bound on $N^*$ gives that if $u_0 = Cq^{1/p}(\sum_{i \leq q} a_i^{p/(p-1)})^{(p-1)/p}$, then

$$\sum_{i \leq q} N^*\left(\frac{qa_i}{u_0}\right) \leq \left(\frac{Cq}{u_0}\right)^{p/(p-1)} \sum_{i \leq q} a_i^{p/(p-1)} \leq q,$$

which yields

$$t_0 \leq Cq^{1/p}\left(\sum_{i \leq q} a_i^{p/(p-1)}\right)^{(p-1)/p}.$$

Moreover, if $i_0$ is the biggest integer in $\{1, \ldots, q+1\}$ such that $qa_{i_0-1}/t_0 \geq 2$, then for all $i \leq i_0 - 1$, $qa_i/t_0 \geq 2$, in which case we can use the lower bound of $N^*$ and for all $i \geq i_0$, $a_i < 2t_0/q$. By definition of $t_0$, we get

$$q \geq \sum_{i \leq q} N^*\left(\frac{qa_i}{t_0}\right) \geq \sum_{i \leq i_0 - 1} \left(\frac{Cqa_i}{t_0}\right)^{p/(p-1)},$$

which shows that $t_0 \geq cq^{1/p}(\sum_{i \leq i_0-1} a_i^{p/(p-1)})^{(p-1)/p}$. It is now clear that

$$q^{1/p}\left(\sum_{i \leq q} a_i^{p/(p-1)}\right)^{(p-1)/p}$$
$$\leq q^{1/p}\left(\sum_{i \leq i_0-1} a_i^{p/(p-1)}\right)^{(p-1)/p} + q^{1/p}\left(\sum_{i \geq i_0} a_i^{p/(p-1)}\right)^{(p-1)/p}$$
$$\leq \frac{t_0}{c} + \left(\frac{q - i_0 + 1}{q}\right)^{(p-1)/p} 2t_0 \leq \left(2 + \frac{1}{c}\right)t_0.$$

Now we establish inequalities (3). To prove the bounds on $\theta_p$, note that since $|g|$ has uniformly bounded density in $p$, there is an absolute constant $c > 0$ such that for every $s > 0$, $P(|g| \geq s) \geq 1 - cs$. If $s = c^{-1}(1 - e^{-1})$, then $P(|g| \geq s) \geq P(|g| \geq 1/\theta_p)$, which shows that $\theta_p \leq s^{-1} \leq C$. On the other hand, Lemma 8 implies that there is an absolute constant $c'$ for which $P(|g| \geq c') \leq 1/e = P(|g| \geq 1/\theta_p)$, and thus $\theta_p \geq 1/c'$.

Finally, we address the above mentioned bounds on $N^*$. Lemma 8 states that $N$ is convex. In particular, $N$ is bounded from below by its tangent function at zero, that is, $N(s) \geq sN'(0)$. So if $t \leq 1/\theta_p \Gamma(1 + 1/p) = N'(0)$, then

$$0 \leq N^*(t) = \sup_{s>0}(ts - N(s)) \leq \sup_{s>0} s(t - N'(0)) = 0,$$

and the claimed upper bound on $N^*$ is obvious. We may restrict attention to $t \geq 1/\theta_p \Gamma(1 + 1/p)$. Denoting $S = s/\theta_p$, Lemma 8 shows that for every $S \geq 1$,

$$N(s) = N(S\theta_p) \geq S^p + (p-1)\log S + \log[p\Gamma(1+1/p)] \geq S^p.$$

Hence, for every $S \geq 1$,

$$ts - N(s) = tS\theta_p - N(S\theta_p) \leq tS\theta_p - S^p$$

$$\leq \sup_{S>0}\{tS\theta_p - S^p\} = (p-1)\left(\frac{t\theta_p}{p}\right)^{p/(p-1)} \leq (Ct)^{p/(p-1)}.$$

For $0 < S < 1$, that is, $0 < s < \theta_p$, $st - N(s) \leq \theta_p t \leq (Ct)^{p/(p-1)}$ since $t$ is bounded from below, and the upper bound for $N^*$ follows.

The lower bound in Lemma 8 shows that there are absolute constants $c, C > 1$ such that if $S \geq c$, $N(s) = N(S\theta_p) \leq (CS)^p$. Therefore $N^*(t) \geq \sup\{tS\theta_p - (CS)^p; S \geq c\}$ and if $t\theta_p \geq c^{p-1}C^p p$, this supremum is attained at $S = (t\theta_p/pC^p)^{1/(p-1)} \geq c$ so that

$$N^*(t) \geq \left(1 - \frac{1}{p}\right)\left(\frac{t\theta_p}{C}\right)^{p/(p-1)} \frac{1}{p^{1/(p-1)}},$$

and we are done. We may therefore assume that $t\theta_p \leq pc^{p-1}C^p$. By our choice of $\theta_p$, $N(1) = 1$, which implies that for all $t \geq 2$,

$$N^*(t) \geq t - N(1) \geq \frac{t}{2} \geq (\tilde{C}t)^{p/(p-1)},$$

with a new constant $\tilde{C}$. This completes the proof. □

The results of this section may be combined to obtain the following exact expression, up to universal constants: for $a_1 \geq a_2 \geq \cdots \geq a_n \geq 0$,

$$(4) \quad \left(\frac{1}{\mathrm{vol}(B_p^n)} \int_{B_p^n} \left|\sum_{i=1}^n a_i x_i\right|^q dx\right)^{1/q} \sim \frac{q^{1/p}\|(a_i)_{i \leq q}\|_{p'} + \sqrt{q}\|(a_i)_{i>q}\|_2}{(\max\{n,q\})^{1/p}},$$

which virtually allows one to solve any question related to moment estimates on $B_p^n$.

2.2.1. *Khinchine inequalities.* A well-known variant of Khinchine's inequality (see [23]) states that for every $1 \leq p, q < \infty$ and every integer $n$, there are $A(p,q,n), B(p,q,n) > 0$ such that for every $(a_1, \ldots, a_n) \in \mathbb{R}^n$,

$$A(p,q,n)\left(\sum_{i=1}^n a_i^2\right)^{1/2} \leq \left(\frac{1}{\mathrm{vol}(B_p^n)}\int_{B_p^n}\left|\sum_{i=1}^n a_i x_i\right|^q dx\right)^{1/q}$$

$$\leq B(p,q,n)\left(\sum_{i=1}^n a_i^2\right)^{1/2},$$

and we assume that $A(p,q,n), B(p,q,n)$ are the best constants for which the above inequality holds for all $(a_1,\ldots,a_n) \in \mathbb{R}^n$. We determine $A(p,q,n)$ and $B(p,q,n)$, up to absolute multiplicative constants.

THEOREM 7. *For every integer $n$ and for every $1 \leq q < \infty$ and $1 \leq p \leq 2$,*

$$A(p,q,n) \sim \frac{\sqrt{q}}{n^{1/p}} \min\left\{1, \sqrt{\frac{n}{q}}\right\} \quad \text{and} \quad B(p,q,n) \sim \min\left\{1, \left(\frac{q}{n}\right)^{1/p}\right\},$$

*while for $2 < p < \infty$,*

$$A(p,q,n) \sim \min\left\{1, \left(\frac{q}{n}\right)^{1/p}\right\} \quad \text{and} \quad B(p,q,n) \sim \frac{\sqrt{q}}{n^{1/p}} \min\left\{1, \sqrt{\frac{n}{q}}\right\}.$$

This is a consequence of (4) and of the following:

LEMMA 9. *For every $a = (a_1,\ldots,a_n) \in S^{n-1}$, if $1 < p \leq 2$ then*

$$\sqrt{q}\max\left\{1, \left(\frac{q}{n}\right)^{1/p - 1/2}\right\} \leq q^{1/p}\left(\sum_{i \leq q} a_i^{p/(p-1)}\right)^{(p-1)/p} + \sqrt{q}\left(\sum_{i > q} a_i^2\right)^{1/2}$$

$$\leq \sqrt{2} \cdot q^{1/p}.$$

*If $2 < p < \infty$, then*

$$q^{1/p} \leq q^{1/p}\left(\sum_{i \leq q} a_i^{p/(p-1)}\right)^{(p-1)/p} + \sqrt{q}\left(\sum_{i > q} a_i^2\right)^{1/2}$$

$$\leq \sqrt{2q}\min\left\{1, \left(\frac{n}{q}\right)^{1/2 - 1/p}\right\}.$$

*Furthermore, these inequalities are optimal, up to universal constants.*

PROOF. Assume that $1 < p \leq 2$. Since $\sqrt{a} + \sqrt{b} \leq \sqrt{2}\sqrt{a+b}$,

$$q^{1/p}\left(\sum_{i \leq q} a_i^{p/(p-1)}\right)^{(p-1)/p} + \sqrt{q}\left(\sum_{i > q} a_i^2\right)^{1/2}$$

$$\leq q^{1/p}\left(\sum_{i \leq q} a_i^2\right)^{1/2} + \sqrt{q}\left(\sum_{i > q} a_i^2\right)^{1/2} \leq \sqrt{2} \cdot q^{1/p}.$$

Similarly, if $q > n$, then

$$q^{1/p}\left(\sum_{i \leq q} a_i^{p/(p-1)}\right)^{(p-1)/p} + \sqrt{q}\left(\sum_{i > q} a_i^2\right)^{1/2} = q^{1/p}\left(\sum_{i \leq n} a_i^{p/(p-1)}\right)^{(p-1)/p}$$

$$\geq \frac{q^{1/p}}{n^{1/p - 1/2}},$$

and if $q \leq n$,

$$q^{1/p}\left(\sum_{i \leq q} a_i^{p/(p-1)}\right)^{(p-1)/p} + \sqrt{q}\left(\sum_{i > q} a_i^2\right)^{1/2} \geq \sqrt{q}.$$

The fact that these inequalities are best possible up to universal constants follows by considering in each case the vectors $(1, 0, \ldots, 0)$, $(1/\sqrt{n}, \ldots, 1/\sqrt{n})$ or $(1/\sqrt{q}, \ldots, 1/\sqrt{q}, 0, \ldots, 0)$ when $q \leq n$. The proof of the case $p \geq 2$ is equally simple. $\square$

2.2.2. *$\psi_2$-directions.* We start with a few definitions. Let $\alpha \in [1, 2]$ and set $\mu$ to be a probability measure on $\mathbb{R}^n$. For a measurable function $f : \mathbb{R}^n \to \mathbb{R}$, define the following Orlicz norm associated with $\alpha$ and $\mu$ by

$$\|f\|_{\psi_\alpha(\mu)} := \inf\left\{\lambda > 0; \int e^{|f/\lambda|^\alpha}\, d\mu \leq 2\right\}.$$

It is well known that $\|f\|_{\psi_\alpha(\mu)} \sim \sup_{q \geq 1} q^{-1/\alpha}(\int |f|^q\, d\mu)^{1/q}$ (this follows from the Taylor expansion of the exponential). Given a vector $\theta$ in the unit sphere $S^{n-1}$ of $\mathbb{R}^n$, one says that $\theta$ defines a $\psi_\alpha$-direction for $\mu$ with a constant $C > 0$ if the function $f_\theta(x) = \langle x, \theta \rangle$ satisfies

$$\|f_\theta\|_{\psi_\alpha(\mu)} \leq C\left(\int |f_\theta|^2\, d\mu\right)^{1/2}.$$

In other words, the moment of $f_\theta$ of order $q$ is bounded from above by a constant times $Cq^{1/\alpha}$ times the second moment of $f_\theta$.

From now on consider a convex body $K \subset \mathbb{R}^n$, with the center of mass at the origin. Such a body is said to be a $\psi_\alpha$-body with constant $C$ if all directions $\theta \in S^{n-1}$ are $\psi_\alpha$ with a constant $C$, with respect to the uniform probability measure on $K$. It follows from the Brunn–Minkowski inequality that convex bodies are $\psi_1$ with a uniform constant, and any improvement on this estimate would be very useful. Note that the notion of $\psi_2$-bodies is crucial in Bourgain's bound on the isotropy constant [12] of convex bodies. This motivated recent works on the $\psi_2$-directions of convex bodies. In fact, it is not even clear that there exists a universal constant $C$ such that any convex body (of any dimension) admits at least one $\psi_2$-direction with constant $C$. This question of Milman was solved in special cases such as zonoids [27] and unconditional bodies (Bobkov and Nazarov [11] show that the main diagonal is $\psi_2$). Thanks to (4) we are able to study these questions for $B_p^n$.

PROPOSITION 10. *There exists $C > 0$ such that:*

(i) *for every $n \geq 1$ and every $p \geq 2$, $B_p^n$ is a $\psi_2$-body with constant $C$.*
(ii) *for every $n \geq 1$ and every $p \in [1, 2]$, $B_p^n$ is a $\psi_p$-body with constant $C$.*

The first point was actually a consequence of results in [8], where subindependence was also used.

PROOF. Without loss of generality we consider a direction $\theta \in S^{n-1}$ with $\theta_1 \geq \theta_2 \geq \cdots \geq \theta_n \geq 0$. Fix $q \geq 1$. Equation (4) gives, with obvious notation,

$$(5) \quad \frac{(\mathbb{E}_{B_p^n}|\langle X, \theta\rangle|^q)^{1/q}}{(\mathbb{E}_{B_p^n}|\langle X, \theta\rangle|^2)^{1/2}} \sim \left(\frac{n}{\max\{n,q\}}\right)^{1/p} \cdot (q^{1/p}\|(\theta_i)_{i\leq q}\|_{p'} + \sqrt{q}\|(\theta_i)_{i>q}\|_2),$$

where $p' = p/(p-1)$.

The result now follows from obvious estimates. Indeed, since $n/\max\{n,q\} \leq 1$, for $p \geq 2$, Hölder's inequality implies that $\|(\theta_i)_{i\leq q}\|_{p'} \leq \min\{n,q\}^{1/2-1/p} \times \|(\theta_i)_{i\leq q}\|_2 \leq q^{1/2-1/p}$. Hence, the right-hand side in (5) is less than $2\sqrt{q}$. For $p \in [1,2]$, it is evident that $\|(\theta_i)_{i\leq q}\|_{p'} \leq \|(\theta_i)_{i\leq q}\|_2 \leq 1$ and thus the ratio of moments is bounded by a constant times $q^{1/p}$. □

Next, we describe the $\psi_2$-constant on $B_p^n$ of every direction for $1 \leq p \leq 2$.

PROPOSITION 11. *Let $p \in [1,2]$. For any integer $n \geq 1$ and $\theta \in S^{n-1}$, $\theta$ is a $\psi_2$-direction of $B_p^n$ and the best constant for which it is $\psi_2$ is, up to an absolute multiplicative constant, $n^{1/p-1/2}\|\theta\|_{p'}$.*

Observe that from the above result, the direction of the main diagonal is $\psi_2$. For $p = 1$ we recover a result of Bobkov and Nazarov [10]. (Let us note that in that paper, the authors give another moment estimate for $B_1^n$, which can be recovered by our method, and which implies that most directions are $\psi_{2-\epsilon}$. Moreover, Bobkov and Nazarov show that these moment upper estimates for $B_1^n$ can be transferred to isotropic unconditional convex bodies.)

PROOF. Assume, as we may, that $\theta_1 \geq \theta_2 \geq \cdots \geq \theta_n \geq 0$. For $q < n$ the right-hand side of (5) is equal to

$$q^{1/p}\|(\theta_i)_{i\leq q}\|_{p'} + \sqrt{q}\|(\theta_i)_{i>q}\|_2 \leq \sqrt{q}(n^{1/p-1/2}\|(\theta_i)_{i\leq n}\|_{p'} + 1)$$
$$\leq 2\sqrt{q}n^{1/p-1/2}\|(\theta_i)_{i\leq n}\|_{p'},$$

where we used Hölder's inequality in the form $1 = \|(\theta_i)_{i\leq n}\|_2 \leq n^{1/p-1/2} \times \|(\theta_i)_{i\leq n}\|_{p'}$.

If $q \geq n$, the right-hand side of (5) is $n^{1/p}\|\theta\|_{p'} \leq \sqrt{q}n^{1/p-1/2}\|(\theta_i)_{i\leq n}\|_{p'}$. For $q = n$, it is easy to see that the estimate cannot be improved by more than a universal factor. □

**3. Extremal geometric parameters of sections of $B_p^n$, $p > 0$.** In what follows we will denote by $G$ a standard Gaussian vector. If $E \subset \mathbb{R}^n$ is a $k$-dimensional subspace, then $G$ will still stand for a standard Gaussian vector on $E$ (which is well defined due to rotational invariance).

3.1. *Bounds via stochastic ordering.* In this section, we present monotonicity properties for sections of $B_p^n$ as $p > 0$ varies. We follow the approach of Meyer and Pajor [22]. They proved that for a fixed vector subspace of dimension $k$ in $\mathbb{R}^n$, the ratio $\mathrm{Vol}_k(E \cap B_p^n)/\mathrm{Vol}_k(B_p^k)$ is nondecreasing in $p \geq 1$. This was later extended to $p > 0$ and to $\ell_p$-sums of arbitrary spaces of finite dimension (see [6] and the reference therein). We are interested in Gaussian averages of the $\ell_p$-norm on sections. Our results will recover in several ways the latter result on the volume.

We will use the notion of peaked ordering on measures. Given two absolutely continuous measures $\mu$ and $\nu$ on $\mathbb{R}^d$, one says that $\nu$ is more peaked than $\mu$ and writes $\mu \prec \nu$ if for every symmetric bounded convex set $C$,

$$\mu(C) \leq \nu(C).$$

In the following statement, we put together the properties that we need. They follow from more general results by Kanter [19].

PROPOSITION 12. *Let $\mu, \nu$ be probability measures on $\mathbb{R}$, with even densities which are nonincreasing on $[0, \infty)$. If $\mu \prec \nu$, then for every $n \geq 1$ one has $\mu^{\otimes n} \prec \nu^{\otimes n}$.*

The aim of the next two lemmas is to relate Gaussian averages of the $\ell_p$-norm on subspaces to the values of some product measures. Let $E \subset \mathbb{R}^n$ be a subspace with $\dim(E) = k$. We denote by $P_E$ the orthogonal projection from $\mathbb{R}^n$ onto $E$ and let $u_{k+1}, \ldots, u_n$ be an orthonormal basis of $E^\perp$. Set

$$B_\infty(E^\perp) = \left\{ x \in E^\perp; \sup_{i=k+1,\ldots,n} |\langle x, u_i \rangle| \leq \tfrac{1}{2} \right\}$$

and for $\epsilon > 0$,

$$E(\epsilon) = \{x \in \mathbb{R}^n; x - P_E(x) \in \epsilon B_\infty(E^\perp)\}.$$

We denote by $\gamma_n$ the standard Gaussian measure on $\mathbb{R}^n$, and by $\gamma_E$ the standard Gaussian distribution on a vector subspace $E$.

LEMMA 13. *Let $E$ be a $k$-dimensional subspace of $\mathbb{R}^n$ and set $h$ to be a continuous function in $L_1(\mathbb{R}^n, \gamma_n)$, with the following property: there exist $K, \eta > 0$ such that for every $x \in \mathbb{R}^n$ one has $|h(x)| \leq K e^{\|x\|_2^2/(2+\eta)}$. Then*

$$\int_E h(x)\, d\gamma_E(x) = \lim_{\epsilon \to 0} \left(\frac{2\pi}{\epsilon^2}\right)^{(n-k)/2} \int_{E(\epsilon)} h(x)\, d\gamma_n(x).$$

PROOF. Fix some $\epsilon > 0$. In the following we recall the dimension of the variable of integration by writing $\int_E f(a) \, d^k a$ when $\dim(E) = k$,

$$(2\pi)^{n/2} \int_{E(\epsilon)} h(x) \, d\gamma_n(x)$$
$$= \int_{\mathbb{R}^n} e^{-\|x\|_2^2/2} \mathbf{1}_{E(\epsilon)}(x) h(x) \, dx$$
$$= \int_{E \times \epsilon B_\infty(E^\perp)} e^{-\|a\|_2^2/2 - \|b\|_2^2/2} h(a + b) \, d^k a \, d^{n-k} b$$
$$= \epsilon^{n-k} \int_{E \times B_\infty(E^\perp)} e^{-\|a\|_2^2/2} e^{-\epsilon^2 \|c\|_2^2/2} h(a + \epsilon c) \, d^k a \, d^{n-k} c.$$

By continuity and dominated convergence, the latter integral converges when $\epsilon$ goes to zero to

$$\mathrm{vol}_{n-k}(B_\infty(E^\perp)) \int_E e^{-\|a\|_2^2/2} h(a) \, d^k a,$$

which gives the claimed result.  □

Fix $0 < p, \lambda < \infty$, let $\alpha(p, \lambda) = 2 \int_0^\infty e^{-\lambda t^p - t^2} \, dt$ and set $\mu_{p,\lambda}$ to be the probability measure on $\mathbb{R}$ defined by

$$d\mu_{p,\lambda}(t) = e^{-\lambda \alpha(p,\lambda)^p |t|^p - \alpha(p,\lambda)^2 t^2} \, dt.$$

LEMMA 14. *Let $E$ be a $k$-dimensional subspace of $\mathbb{R}^n$ and $0 < p, \lambda < \infty$. Then*

$$\frac{\mathbb{E} e^{(-\lambda/2^{p/2})\|G\|_{E \cap B_p^n}^p}}{\mathbb{E} e^{(-\lambda/2^{p/2})\|G\|_{B_p^k}^p}} = \lim_{\epsilon \to 0} \epsilon^{k-n} \mu_{p,\lambda}^{\otimes n}(E(\epsilon)).$$

PROOF. By Lemma 13,

$$\mathbb{E} e^{(-\lambda/2^{p/2})\|G\|_{E \cap B_p^n}^p}$$
$$= \lim_{\epsilon \to 0} \left(\frac{2\pi}{\epsilon^2}\right)^{(n-k)/2} \int_{E(\epsilon)} e^{(-\lambda\|x\|_p^p/2^{p/2}) - (\|x\|_2^2/2)} \, dx$$
$$= \lim_{\epsilon \to 0} \left(\frac{2\pi}{\epsilon^2}\right)^{(n-k)/2} 2^{n/2} \alpha(p, \lambda)^n$$
$$\quad \times \int_{E(\epsilon/(\sqrt{2}\alpha(p,\lambda)))} e^{-\lambda \alpha(p,\lambda)^p \|x\|_p^p - \alpha(p,\lambda)^2 \|x\|_2^2} \, dx$$
$$= \lim_{\epsilon \to 0} 2^{k/2} \alpha(p, \lambda)^k \left(\frac{2\pi}{\epsilon^2}\right)^{(n-k)/2} \mu_{p,\lambda}^{\otimes n}(E(\epsilon)).$$

Thus, applied to $E = \{x \in \mathbb{R}^n; x_{k+1} = \cdots = x_n = 0\}$ with $u_i = e_i$ for $i > k$, this identity yields

$$\mathbb{E}e^{(-\lambda/2^{p/2})\|G\|_{B_p^k}^p}$$
$$= \lim_{\epsilon \to 0} 2^{k/2} \alpha(p,\lambda)^k \left(\frac{2\pi}{\epsilon^2}\right)^{(n-k)/2} \left(\int_{-\epsilon/2}^{\epsilon/2} e^{-\lambda \alpha(p,\lambda)^p |t|^p - \alpha(p,\lambda)^2 t^2} \, dt\right)^{n-k}$$
$$= 2^{k/2} \alpha(p,\lambda)^k (2\pi)^{(n-k)/2},$$

from which the required result follows. □

In the forthcoming lemmas and propositions, we look for comparison results in the sense of the peaked ordering, between measures of the form $\mu_{p,\lambda}$. We start with useful facts about the constants $\alpha(p,\lambda)$ which appear in the definition of $\mu_{p,\lambda}$.

LEMMA 15. *Let $\lambda > 0$ and $0 < p < q < \infty$. Then $\alpha(p, \frac{\lambda}{\Gamma((p+1)/2)}) < \alpha(q, \frac{\lambda}{\Gamma((q+1)/2)})$.*

PROOF. By its definition,

$$\alpha(p,\lambda) = 2\int_0^\infty e^{-\lambda t^p - t^2} \, dt = \sqrt{\pi} \cdot \mathbb{E}\exp\left(-\frac{\lambda|g|^p}{2^{p/2}}\right),$$

where $g$ is a standard Gaussian random variable. Recall that

$$\mathbb{E}|g|^p = \frac{2^{p/2}}{\sqrt{\pi}} \Gamma\left(\frac{p+1}{2}\right),$$

and thus

$$\alpha\left(p, \frac{\lambda}{\Gamma((p+1)/2)}\right) = \sqrt{\pi} \cdot \mathbb{E}\exp\left(-\frac{\lambda|g|^p}{\sqrt{\pi} \cdot \mathbb{E}|g|^p}\right). \qquad \square$$

Therefore, Lemma 15 follows from the following result:

LEMMA 16. *Fix $0 < p < q < \infty$ and let $X$ be a nonnegative random variable with $\mathbb{E}X^q < \infty$. Then for every convex function $f: [0,\infty] \to [0,\infty)$,*

$$\mathbb{E}f\left(\frac{X^p}{\mathbb{E}X^p}\right) \leq \mathbb{E}f\left(\frac{X^q}{\mathbb{E}X^q}\right).$$

PROOF. Let $t_0$ be defined by $t_0^{1/p}(\mathbb{E}X^p)^{1/p} = t_0^{1/q}(\mathbb{E}X^q)^{1/q}$. Clearly,

$$(6) \quad \mathbb{E}\left[f\left(\frac{X^p}{\mathbb{E}X^p}\right) - f\left(\frac{X^q}{\mathbb{E}X^q}\right)\right] = \int_0^{t_0} f'(t)h(t) \, dt + \int_{t_0}^\infty f'(t)h(t) \, dt,$$

where $h(t) = P(X^p \geq t\mathbb{E}X^p) - P(X^q \geq t\mathbb{E}X^q)$. Since $h \geq 0$ on $[0, t_0]$ and $h \leq 0$ on $[t_0, \infty)$ and $\int_0^\infty h(t)\,dt = 0$, then

$$\int_0^{t_0} f'(t)h(t)\,dt + \int_{t_0}^\infty f'(t)h(t)\,dt$$

$$= \int_0^{t_0}[f'(t) - f'(t_0)]h(t)\,dt + \int_{t_0}^\infty [f'(t) - f'(t_0)]h(t)\,dt \leq 0,$$

where we have used the fact that $f'$ is nondecreasing. Combined with (6), this completes the proof. □

PROPOSITION 17. *Let $0 < p < q$ and $\lambda_1, \lambda_2 > 0$. Then:*

(a) *If $q \geq 2$ and $\alpha(p, \lambda_1) > \alpha(q, \lambda_2)$, then $\mu_{p,\lambda_1} \prec \mu_{q,\lambda_2}$.*
(b) *If $q < 2$ and $\alpha(p, \lambda_1) < \alpha(q, \lambda_2)$, then $\mu_{p,\lambda_1} \prec \mu_{q,\lambda_2}$.*
(c) *If $p < 2$ and $q \geq 2$, then without any restriction on $\lambda_1$ and $\lambda_2$, $\mu_{p,\lambda_1} \prec \mu_{q,\lambda_2}$.*
(d) *If $0 < p < 2$ and $\lambda_1 < \lambda_2$, then $\mu_{p,\lambda_2} \prec \mu_{p,\lambda_1}$.*
(e) *If $p > 2$ and $\lambda_1 < \lambda_2$, then $\mu_{p,\lambda_1} \prec \mu_{p,\lambda_2}$.*

PROOF. Define $h : [0, \infty) \to \mathbb{R}$ by

$$h(a) = \int_0^a [e^{-\lambda_1\alpha(p,\lambda_1)^p t^p - \alpha(p,\lambda_1)^2 t^2} - e^{-\lambda_2\alpha(q,\lambda_2)^q t^q - \alpha(q,\lambda_2)^2 t^2}]\,dt.$$

In order to prove that $\mu_{p,\lambda_1} \prec \mu_{q,\lambda_2}$ one has to show that $h(a) \leq 0$ for all $a \geq 0$. Note that $h(0) = \lim_{x \to \infty} h(x) = 0$, and if

$$\psi(t) = -\lambda_1\alpha(p,\lambda_1)^p t^{p-2} - \alpha(p,\lambda_1)^2 + \lambda_2\alpha(q,\lambda_2)^q t^{q-2} + \alpha(q,\lambda_2)^2,$$

then $\text{sign}(h') = \text{sign}(\psi)$.

In case (a), $\lim_{a \to 0} \psi(a) < 0$ and $\lim_{a \to \infty} \psi(a) > 0$. Hence $h' < 0$ in a neighborhood of 0 and $h'(a) > 0$ for $a$ large enough. If there were some $a_0 > 0$ such that $h(a_0) > 0$, then it would follow that $h'$ must have at least three zeros. Thus $\psi$ would also have three zeros, implying that $\psi'$ has at least two zeros. This is impossible since

$$\psi'(t) = -\lambda_1(p-2)\alpha(p,\lambda_1)^p t^{p-3} + \lambda_2(q-2)\alpha(q,\lambda_2)^q t^{q-3}$$

clearly has at most one zero.

Cases (b) and (c) are just as simple. To prove case (d) one must show that the function

$$\psi(t) = (\lambda_2\alpha(p,\lambda_2)^p - \lambda_1\alpha(p,\lambda_1)^p)t^{p-2} + \alpha(p,\lambda_2)^2 - \alpha(p,\lambda_1)^2$$

is first positive and then negative. Since it changes signs only once, it is enough to check this at zero and infinity. Observe that

$$\alpha(p,\lambda) = 2\int_0^\infty e^{-\lambda t^p - t^2}\,dt \quad \text{and} \quad \lambda^{1/p}\alpha(p,\lambda) = 2\int_0^\infty e^{-t^p - t^2/\lambda^{2/p}}\,dt,$$

so that $\alpha(p, \lambda)$ is decreasing in $\lambda$ and $\lambda\alpha(p,\lambda)^p$ is increasing in $\lambda$. Since $p < 2$, then $\lim_{x \to 0} \psi(x) = +\infty$ and $\lim_{x \to \infty} \psi(x) < 0$. The proof of the last case is almost identical. $\square$

PROPOSITION 18. *Let $E$ be a $k$-dimensional subspace of $\mathbb{R}^n$ and set $\lambda > 0$. For $p > 0$, let*

$$F(p) = \frac{\mathbb{E}\exp[-(\lambda\|G\|_{E\cap B_p^n}^p)/(2^{p/2}\Gamma((p+1)/2))]}{\mathbb{E}\exp[-(\lambda\|G\|_{B_p^k}^p)/(2^{p/2}\Gamma((p+1)/2))]}.$$

*Then $F$ is nondecreasing on $(0,2]$. Moreover, for $p \geq 2$ one has $F(p) \geq F(2) = 1$.*

PROOF. Let $r < 2$, fix some $\lambda > 0$, let $p > r$ and define

$$\lambda_1 = \frac{\lambda}{\Gamma((r+1)/2)} \quad \text{and} \quad \lambda_2 = \frac{\lambda}{\Gamma((p+1)/2)}.$$

By Lemma 15 and cases (b) and (c) of Proposition 17, $\mu_{r,\lambda_1} \prec \mu_{p,\lambda_2}$. Tensorizing and applying Proposition 12, it follows that $\mu_{r,\lambda_1}^{\otimes n} \prec \mu_{p,\lambda_2}^{\otimes n}$. In particular, for every $\epsilon > 0$,

$$\mu_{r,\lambda_1}^{\otimes n}(E(\epsilon)) \leq \mu_{p,\lambda_2}^{\otimes n}(E(\epsilon)).$$

By Lemma 14,

(7)
$$\frac{\mathbb{E}\exp[-(\lambda\|G\|_{E\cap B_r^n}^r)/(2^{r/2}\Gamma((r+1)/2))]}{\mathbb{E}\exp[-(\lambda\|G\|_{B_r^k}^r)/(2^{r/2}\Gamma((r+1)/2))]}$$
$$\leq \frac{\mathbb{E}\exp[-(\lambda\|G\|_{E\cap B_p^n}^p)/(2^{p/2}\Gamma((p+1)/2))]}{\mathbb{E}\exp[-(\lambda\|G\|_{B_p^k}^p)/(2^{p/2}\Gamma((p+1)/2))]},$$

hence $F(r) \leq F(p)$ holds when $r < 2$ and $r < p$. $\square$

THEOREM 8. *Let $E$ be a $k$-dimensional subspace of $\mathbb{R}^n$. Then the function*

$$p \mapsto \frac{\mathbb{E}\|G\|_{E\cap B_p^n}^p}{\mathbb{E}\|G\|_{B_p^k}^p}$$

*is nonincreasing in $p > 0$.*

PROOF. Assume that $p < q \leq 2$. Both sides of (7) equal 1 for $\lambda = 0$, so the same inequality must hold between the derivatives at 0 of both sides;

that is,

$$-\frac{\mathbb{E}\|G\|^p_{E \cap B^n_p}}{2^{p/2}\Gamma((p+1)/2)} + \frac{\mathbb{E}\|G\|^p_{B^k_p}}{2^{p/2}\Gamma((p+1)/2)}$$
$$\leq -\frac{\mathbb{E}\|G\|^q_{E \cap B^n_q}}{2^{q/2}\Gamma((q+1)/2)} + \frac{\mathbb{E}\|G\|^q_{B^k_q}}{2^{q/2}\Gamma((q+1)/2)}.$$

Note that

$$\mathbb{E}\|G\|^p_{B^k_p} = \mathbb{E}\sum_{i=1}^k |g_i|^p = \frac{2k}{\sqrt{2\pi}} \int_0^\infty x^p e^{-x^2/2}\,dx = \frac{2^{p/2+1}}{\sqrt{\pi}} k\Gamma\left(\frac{p+1}{2}\right).$$

Hence, the above inequality translates to

$$-\frac{2k}{\sqrt{\pi}} \cdot \frac{\mathbb{E}\|G\|^p_{E \cap B^n_p}}{\mathbb{E}\|G\|^p_{B^k_p}} + \frac{2k}{\sqrt{\pi}} \leq -\frac{2k}{\sqrt{\pi}} \cdot \frac{\mathbb{E}\|G\|^q_{E \cap B^n_q}}{\mathbb{E}\|G\|^q_{B^k_q}} + \frac{2k}{\sqrt{\pi}},$$

so that

$$\frac{\mathbb{E}\|G\|^p_{E \cap B^n_p}}{\mathbb{E}\|G\|^p_{B^k_p}} \geq \frac{\mathbb{E}\|G\|^q_{E \cap B^n_q}}{\mathbb{E}\|G\|^q_{B^k_q}}.$$

It remains to deal with the case $2 \leq p < q$, which is slightly more complicated because the last proposition does not give much in this case for a fixed value of the parameter $\lambda$. However, something remains true when $\lambda$ tends to zero, and thus one can pass to the limit.

Indeed, fix two numbers $c_p, c_q > 0$ such that

$$c_p < \frac{1}{\Gamma((p+1)/2)} \quad \text{and} \quad c_q > \frac{1}{\Gamma((q+1)/2)},$$

and for every $\lambda > 0$ define

$$f(\lambda) = \alpha(p, c_p\lambda) - \alpha(q, c_q\lambda) = 2\int_0^\infty e^{-c_p\lambda t^p - t^2}\,dt - 2\int_0^\infty e^{-c_q\lambda t^q - t^2}\,dt.$$

Then

$$f'(0) = -2c_p \int_0^\infty t^p e^{-t^2}\,dt + 2c_q \int_0^\infty t^q e^{-t^2}\,dt$$
$$= 2 \cdot \left[c_q\Gamma\left(\frac{q+1}{2}\right) - c_p\Gamma\left(\frac{p+1}{2}\right)\right] > 0.$$

Since $f(0) = 0$, it follows that there is some $\delta = \delta_{p,q} > 0$ such that for every $0 < \lambda < \delta$, $f(\lambda) > 0$, that is, $\alpha(p, c_p\lambda) > \alpha(q, c_q\lambda)$. Part (a) of Proposition 17 now implies that $\mu_{p,c_p\lambda} \prec \mu_{q,c_q\lambda}$. As before, tensorization and an application of Lemma 14 give that for every $\lambda < \delta$,

$$\frac{\mathbb{E}\exp(-\lambda c_p\|G\|^p_{E\cap B_p^n}/2^{p/2})}{\mathbb{E}\exp(-\lambda c_p\|G\|^p_{B_p^k}/2^{p/2})} \leq \frac{\mathbb{E}\exp(-\lambda c_q\|G\|^q_{E\cap B_q^n}/2^{q/2})}{\mathbb{E}\exp(-\lambda c_q\|G\|^q_{B_q^k}/2^{q/2})},$$

and the required inequality follows by taking derivatives at 0 and letting $c_p$ and $c_q$ tend to $1/\Gamma((p+1)/2)$ and $1/\Gamma((q+1)/2)$, respectively. $\square$

REMARK. Assume that $0 < p < 2$. By Proposition 17, for every $\lambda > 0$, $\mu_{p,\lambda} \leq \mu_{2,\lambda} = \bar{\gamma}$, where $\bar{\gamma}$ has density $e^{-\pi x^2}$ on $\mathbb{R}$. Hence, by rotation invariance of this Gaussian density, one has that for every $\lambda > 0$,

$$\mathbb{E}e^{-\lambda\|G\|^p_{E\cap B_p^n}} \leq \mathbb{E}e^{-\lambda\|G\|^p_{B_p^k}}.$$

Thus, for any (reasonable) measure $\tau$ on $[0,\infty)$,

$$\mathbb{E}\int_0^\infty e^{-\lambda\|G\|^p_{E\cap B_p^n}}\, d\tau(\lambda) \leq \mathbb{E}\int_0^\infty e^{-\lambda\|G\|^p_{B_p^k}}\, d\tau(\lambda),$$

which by Bernstein's theorem (see, e.g., [36]) implies that for every $f:[0,\infty)\to\mathbb{R}$ which is completely monotonic,

$$\mathbb{E}f(\|G\|^p_{E\cap B_p^n}) \leq \mathbb{E}f(\|G\|^p_{B_p^k}),$$

provided these expectations are finite.

Two particular cases which should be singled out are $f(t) = e^{-\lambda t^\theta}$ for $0 < \theta \leq 1$ and $\lambda > 0$, and $f(t) = t^{-\eta}$ for $\eta > 0$. The first case implies that for every $\lambda > 0$,

$$\mathbb{E}e^{-\lambda\|G\|^{\theta p}_{E\cap B_p^n}} \leq \mathbb{E}e^{-\lambda\|G\|^{\theta p}_{B_p^k}},$$

which by differentiation at 0 yields

$$\mathbb{E}\|G\|^{\theta p}_{E\cap B_p^n} \geq \mathbb{E}\|G\|^{\theta p}_{B_p^k}.$$

From the second case it is evident that for $0 < \alpha < k$,

$$\mathbb{E}\|G\|^{-\alpha}_{E\cap B_p^n} \leq \mathbb{E}\|G\|^{-\alpha}_{B_p^k}.$$

The condition $\alpha < k$ is imposed to ensure that these expectations would be finite.

When $2 < p < \infty$, $\bar{\gamma} \prec \mu_{p,\lambda}$, and all the above inequalities are reversed. Summarizing, we obtain

COROLLARY 19. *Let $E$ be a $k$-dimensional subspace of $\mathbb{R}^n$. Then for $0 < p < 2$ and every $0 < \alpha < k$ and $0 < \beta \leq p$,*

$$\mathbb{E}\|G\|_{E \cap B_p^n}^{-\alpha} \leq \mathbb{E}\|G\|_{B_p^k}^{-\alpha} \quad \text{and} \quad \mathbb{E}\|G\|_{E \cap B_p^n}^{\beta} \geq \mathbb{E}\|G\|_{B_p^k}^{\beta}.$$

*If $2 < p < \infty$, then for every $0 < \alpha < k$ and $0 < \beta \leq p$,*

$$\mathbb{E}\|G\|_{E \cap B_p^n}^{-\alpha} \geq \mathbb{E}\|G\|_{B_p^k}^{-\alpha} \quad \text{and} \quad \mathbb{E}\|G\|_{E \cap B_p^n}^{\beta} \leq \mathbb{E}\|G\|_{B_p^k}^{\beta}.$$

The following proposition is a corollary of parts (d) and (e) in Proposition 17.

PROPOSITION 20. *Let $E$ be a $k$-dimensional subspace of $\mathbb{R}^n$. Then the function*

$$\lambda \geq 0 \mapsto r_p(\lambda) := \frac{\mathbb{E} e^{-\lambda \|G\|_{E \cap B_p^n}^p}}{\mathbb{E} e^{-\lambda \|G\|_{B_p^k}^p}}$$

*is nonincreasing when $p \leq 2$ and nondecreasing when $p \geq 2$.*

REMARK. Since $r_p(0) = 1$, we have an alternative proof to Corollary 19. Additionally, the limit of $r_p(\lambda)$ when $\lambda$ tends to infinity is

$$\frac{\int_E e^{-\|x\|_{E \cap B_p^n}^p} \, dx}{\int_{\mathbb{R}^k} e^{-\|x\|_{B_p^k}^p} \, dx} = \frac{\mathrm{vol}_k(E \cap B_p^n)}{\mathrm{vol}_k(B_p^k)}.$$

The above equality can be proved by polar integration. The comparison between $r_p(0)$ and $r_p(+\infty)$ yields an alternative proof of the Meyer–Pajor theorem [22] which uses a different interpolation between $\exp(-t^2)$ and $\exp(-|t|^p)$.

3.2. *Bounds via convolution inequalities.* In this section we derive upper bounds on the Laplace transform of $\|G\|_{E \cap B_p^n}^p$ for $p > 2$. The main tool is Ball's version of the Brascamp–Lieb inequality [3, 13]. We follow the method of [3] where the main focus was on the volume of sections.

Let $E$ be a $k$-dimensional subspace of $\mathbb{R}^n$ and let $P$ be the orthogonal projection onto $E$. The canonical basis of $\mathbb{R}^n$ provides a decomposition of the identity map as $\sum_{i=1}^n e_i \otimes e_i = \mathrm{Id}_n$, where $(v \otimes v)(x) = \langle x, v \rangle v$. Projecting this relation onto $E$ yields a decomposition of the identity on $E$

$$\sum_{i=1}^n Pe_i \otimes Pe_i = \mathrm{Id}_E.$$

Setting $c_i = |Pe_i|^2$ and $u_i = Pe_i/|Pe_i|$ (or any unit vector if the norm of $Pe_i$ is 0), this rewrites as $\sum_{i=1}^n c_i u_i \otimes u_i = \mathrm{Id}_E$. Let $\lambda > 0$, and note that for any $x \in E$ the $i$th coordinate in the canonical basis is $x_i = \langle x, e_i \rangle = \langle Px, e_i \rangle = \langle x, Pe_i \rangle = \sqrt{c_i} \langle x, u_i \rangle$. Hence,

$$\int_E e^{-\lambda \|x\|_p^p - \|x\|_2^2/2}\, dx = \int_E \prod_{i=1}^n e^{-\lambda |x_i|^p - |x_i|^2/2}\, dx$$

$$= \int_E \prod_{i=1}^n e^{-\lambda |\sqrt{c_i}\langle x, u_i\rangle|^p - c_i \langle x, u_i\rangle^2/2}\, dx$$

$$= \int_E \prod_{i=1}^n (e^{-\lambda c_i^{p/2-1} |\langle x, u_i\rangle|^p - \langle x, u_i\rangle^2/2})^{c_i}\, dx$$

$$\leq \prod_{i=1}^n \left( \int_{\mathbb{R}} e^{-\lambda c_i^{p/2-1}|t|^p - t^2/2}\, dt \right)^{c_i}$$

$$= \exp\left[ \sum_{i=1}^n c_i \log \psi\left( \frac{1}{\sqrt{c_i}} \right) \right],$$

where we have set $\psi(s) = 2\int_0^\infty e^{-\lambda s^{2-p} t^p - t^2/2}\, dt$. First, observe that for $p > 2$ the function defined on $(0,\infty) \times [0,\infty)$ by $(s,t) \to -\lambda s^{2-p} t^p - t^2/2$ is concave. [Indeed, on this set the function $(s,t) \to s^{2-p} t^p$ is convex, as follows from a direct calculation of its Hessian matrix.] Therefore, by a well-known result of Borell, Prékopa and Rinott (see, e.g., [30]), $\log \psi(s)$ is a concave function of $s > 0$ [because it is the integral in $t$ of a log-concave function of $(s,t)$]. Lemma 22 below ensures that the map

$$s > 0 \mapsto s \log \psi\left( \frac{1}{\sqrt{s}} \right)$$

is concave. This property can be combined with the relation $\sum_{i=1}^n c_i = k$ (which follows by taking traces in the decomposition of the identity). It yields that for $p \geq 2$,

$$\int_E e^{-\lambda \|x\|_p^p - \|x\|_2^2/2}\, dx \leq \left( \int_{\mathbb{R}} e^{-\lambda (\sqrt{k/n})^{p-2} |t|^p - t^2/2}\, dt \right)^k.$$

Returning to our previous setting, it implies that for every $\lambda > 0$,

$$\mathbb{E} e^{-\lambda \|G\|_{E \cap B_p^n}^p} \leq \mathbb{E} e^{-\lambda(\sqrt{k/n})^{p-2} \|G\|_{B_p^k}^p}.$$

Integrating this inequality against positive measures on $[0, \infty)$ and applying Bernstein's theorem [36], it follows that for every completely monotonic function $f : [0, \infty) \to [0, \infty)$,

$$\mathbb{E} f(\lambda \|G\|_{E \cap B_p^n}^p) \leq \mathbb{E} f(\lambda (\sqrt{k/n})^{p-2} \|G\|_{B_p^k}^p).$$

In particular, the following corollary is evident.

COROLLARY 21. *For any $p \geq 2$, every $0 \leq \theta \leq 1$ and every $\lambda \geq 0$,*

$$\mathbb{E} e^{-\lambda \|G\|_{E \cap B_p^n}^{\theta p}} \leq \mathbb{E} e^{-\lambda (\sqrt{k/n})^{\theta(p-2)} \|G\|_{B_p^k}^{\theta p}}.$$

*In particular, by differentiation at 0 it follows that for every $0 \leq \beta \leq p$,*

$$\mathbb{E}\|G\|_{E \cap B_p^n}^{\beta} \geq \left(\frac{k}{n}\right)^{\beta(1/2 - 1/p)} \mathbb{E}\|G\|_{B_p^k}^{\beta}.$$

*Also, for every $0 \leq \alpha < k$,*

$$\mathbb{E}\|G\|_{E \cap B_p^n}^{-\alpha} \leq \left(\frac{n}{k}\right)^{\alpha(1/2 - 1/p)} \mathbb{E}\|G\|_{B_p^k}^{-\alpha}.$$

REMARK. Assume that $k$ divides $n$, and write $n = mk$. Consider the subspace $F \subset \mathbb{R}^n$ which is the "main diagonal" with respect to the decomposition $\mathbb{R}^n = \mathbb{R}^k \times \cdots \times \mathbb{R}^k$ [i.e., $F = \{(x_1, \ldots, x_m); x_i \in \mathbb{R}^k, x_1 = \cdots = x_m\}$]. Then

$$\mathbb{E}\|G\|_{F \cap B_p^n}^p = m\left(\frac{1}{\sqrt{m}}\right)^p \mathbb{E}\|G\|_{B_p^k}^p = \left(\frac{k}{n}\right)^{p/2 - 1} \mathbb{E}\|G\|_{B_p^k}^p,$$

which shows that when $k$ divides $n$, the case $\beta = p$ in Corollary 21 is optimal.

LEMMA 22. *Let $c:[0, \infty) \to [0, \infty)$ be a nondecreasing concave function. Then the function $f(t) := tc(\frac{1}{\sqrt{t}})$, defined for $t > 0$, is concave.*

PROOF. We may assume that $c$ is twice continuously differentiable. Clearly,

$$f'(t) = c\left(\frac{1}{\sqrt{t}}\right) - \frac{1}{2\sqrt{t}} c'\left(\frac{1}{\sqrt{t}}\right),$$

which is nonincreasing provided the function $g(u) = c(u) - \frac{u}{2} c'(u)$ is nondecreasing on $[0, \infty)$. Now, $g'(u) = \frac{c'(u)}{2} - \frac{u}{2} c''(u)$ is nonnegative by our assumptions on $c$. $\square$

3.3. *Gaussian measures of sections of the cube.* In view of the previous results, one is tempted to conjecture that the following distributional inequality holds for Gaussian measures of sections of dilates of the $\ell_p^n$-ball, that is, for every $k$-dimensional subspace $E$ and every $r > 0$, $\gamma_k(rB_p^k) \leq \gamma_E(E \cap rB_p^n)$ if $p \geq 2$ and the reverse inequality for $p \leq 2$. If such a statement were true, some of the previous results would follow by integration. Unfortunately, it seems that the known techniques are insufficient for this purpose. The product structure of the cube will, however, allow us to prove this conjecture for $p = \infty$.

By Lemma 13, for every $k$-dimensional subspace $E \subset \mathbb{R}^n$ and $r > 0$,

$$\gamma_E(E \cap rB_\infty^n) = \lim_{\epsilon \to 0} \left(\frac{\pi}{2\epsilon^2}\right)^{(n-k)/2} \frac{1}{(2\pi)^{n/2}} \int_{E(\epsilon)} \prod_{i=1}^n e^{-x_i^2/2} \mathbf{1}_{[-r,r]}(x_i) \, dx.$$

Let $\theta(r) = \theta$ be such that

$$\int_{-r/\theta}^{r/\theta} e^{-\theta^2 t^2/2} \, dt = 1,$$

that is,

$$\theta(r) = \int_{-r}^{r} e^{-t^2/2} \, dt.$$

Clearly $\theta$ is increasing and the function $r \mapsto \frac{\theta(r)}{r}$ is decreasing.

Denote by $\rho_r$ the probability measure on $\mathbb{R}$ defined by

$$d\rho_r(t) = e^{-\theta(r)^2 t^2/2} \mathbf{1}_{[-r/\theta(r), r/\theta(r)]}(t) \, dt.$$

Thus,

$$\gamma_E(E \cap rB_\infty^n)$$
$$= \lim_{\epsilon \to 0} \left(\frac{\pi}{2\epsilon^2}\right)^{(n-k)/2} \frac{\theta(r)^n}{(2\pi)^{n/2}} \int_{E(\epsilon/\theta(r))} \prod_{i=1}^n e^{-\theta(r)^2 y_i^2/2} \mathbf{1}_{[-r,r]}(\theta(r) y_i) \, dy$$
$$= \lim_{\epsilon \to 0} \left(\frac{2}{\pi}\right)^{k/2} \cdot \frac{\theta(r)^k}{2^n \epsilon^{n-k}} \cdot \rho_r^{\otimes n}(E(\epsilon)).$$

Observe that

$$\gamma_k(rB_\infty^k) = \frac{1}{(2\pi)^{k/2}} \left(\int_{-r}^{r} e^{-t^2/2} \, dt\right)^k = \frac{\theta(r)^k}{(2\pi)^{k/2}},$$

hence

(8) $$\frac{\gamma_E(E \cap rB_\infty^n)}{\gamma_k(rB_\infty^k)} = \lim_{\epsilon \to 0} \frac{1}{(2\epsilon)^{n-k}} \cdot \rho_r^{\otimes n}(E(\epsilon)).$$

LEMMA 23. *For every $r > s > 0$, $\rho_r \prec \rho_s$.*

PROOF. As usual, define $h: [0, \infty) \to \mathbb{R}$ by

$$h(a) = \int_0^a [e^{-\theta(r)^2 t^2/2} \mathbf{1}_{[-r/\theta(r), r/\theta(r)]}(t) - e^{-\theta(s)^2 t^2/2} \mathbf{1}_{[-s/\theta(s), s/\theta(s)]}(t)] \, dt,$$

and our goal is to show that $h(a) \leq 0$ for all $a \geq 0$. The above mentioned properties of $\theta$ yield $\frac{r}{\theta(r)} \geq \frac{s}{\theta(s)}$, so that $h(a) = 0$ for $a \geq \frac{r}{\theta(r)}$. Moreover, for $\frac{s}{\theta(s)} \leq a \leq \frac{r}{\theta(r)}$, $h(a) = \rho_r([0,a]) - 1 \leq 0$. Finally, for $0 \leq a \leq \frac{s}{\theta(s)}$,

$$h(a) = \int_0^a [e^{-\theta(r)^2 t^2/2} - e^{-\theta(s)^2 t^2/2}]\, dt \leq 0,$$

since $\theta(r) \geq \theta(s)$.  □

By (8), tensorizing the above lemma yields:

THEOREM 9.  *For every $k$-dimensional subspace $E \subset \mathbb{R}^n$ the function*

$$r \mapsto \frac{\gamma_E(E \cap rB_\infty^n)}{\gamma_k(rB_\infty^k)}, \qquad r > 0,$$

*is nonincreasing. In particular, by passing to the limit $r \to \infty$ it follows that for every $r > 0$,*

$$\gamma_E(E \cap rB_\infty^n) \geq \gamma_k(rB_\infty^k).$$

By arguments analogous to those in Section 3.2 one can also obtain the following upper bound on the Gaussian measure of sections of dilates of the cube, which is a Gaussian analog of Ball's slicing theorem in [3]. As noted in Section 3.2, these bounds are optimal when $k$ divides $n$.

THEOREM 10.  *For every $k$-dimensional subspace $E \subset \mathbb{R}^n$ and every $r > 0$,*

$$\gamma_E(E \cap rB_\infty^n) \leq \gamma_k\left(r\sqrt{\frac{n}{k}} B_\infty^k\right).$$

3.4. *An application: a remark on the Komlós conjecture.* In this section we apply the results of the previous section to prove the following proposition, which was stated in the Introduction:

PROPOSITION 24.  *There is an absolute constant $C > 0$ such that for every integer $m > 0$ and every $x_1, \ldots, x_m \in \ell_\infty$, if we denote by $d$ the dimension of the linear span of $x_1, \ldots, x_m$, then there are signs $\varepsilon_1, \ldots, \varepsilon_m \in \{-1, 1\}$ such that*

$$\left\| \sum_{i=1}^m \varepsilon_i x_i \right\|_\infty \leq C\sqrt{\log d} \cdot \max_{1 \leq i \leq m} \|x_i\|_2 \leq C\sqrt{\log m} \cdot \max_{1 \leq i \leq m} \|x_i\|_2.$$

PROOF. We may assume that $x_1,\ldots,x_m \in \ell_2$, in which case we may write $x_i = y_i + z_i$, where $y_i \in \ell_\infty^N$ for some (large) $N$, and $\|z_i\|_\infty \leq 1/m$. Denote $E = \text{span}\{y_1,\ldots,y_m\}$ and let $d'$ be the dimension of $E$. There is a constant $c > 0$ such that for $r = c\sqrt{\log d'} \leq c\sqrt{\log d}$, $\gamma_{d'}(rB_\infty^{d'}) \geq \frac{1}{2}$. By Theorem 9, if we set $K = E \cap rB_\infty^N$, then $\gamma_E(K) \geq \frac{1}{2}$. By Banaszczyk's theorem [5], there are signs $\varepsilon_1,\ldots,\varepsilon_m \in \{-1,1\}$ such that $\sum_{i=1}^m \varepsilon_i y_i \in cK$, where $c$ is an absolute constant. Hence

$$\left\|\sum_{i=1}^m \varepsilon_i x_i\right\|_\infty \leq \left\|\sum_{i=1}^m \varepsilon_i y_i\right\|_\infty + \sum_{i=1}^m \|z_i\|_\infty \leq (c+1)\sqrt{\log d}. \qquad \square$$

It is equally simple to deduce the following $\ell_p$-version of this result for $p > 2$:

PROPOSITION 25. *There is an absolute constant $C > 0$ such that for every $2 \leq p < \infty$, every integer $m > 0$ and every $x_1,\ldots,x_m \in \ell_p$, if we denote by $d$ the dimension of the linear span of $x_1,\ldots,x_m$, then there are signs $\varepsilon_1,\ldots,\varepsilon_m \in \{-1,1\}$ such that*

$$\left\|\sum_{i=1}^m \varepsilon_i x_i\right\|_p \leq C\sqrt{p} \cdot d^{1/p} \cdot \max_{1 \leq i \leq m} \|x_i\|_2 \leq C\sqrt{p} \cdot m^{1/p} \cdot \max_{1 \leq i \leq m} \|x_i\|_2.$$

PROOF. As before, we may assume that $x_1,\ldots,x_m \in \ell_\infty^N$ for some large $N$. By Corollary 19, if we set $E = \text{span}\{x_1,\ldots,x_m\}$, then

$$\mathbb{E}\|G\|_{E \cap B_p^N}^p \leq \mathbb{E}\|G\|_{B_p^d}^p = d\mathbb{E}|g_1|^p = O(dp^{p/2}).$$

Hence, for every $r > 0$,

$$\gamma_E(E \cap rB_p^N) = 1 - P(\|G\|_{E \cap B_p^N}^p \geq r^p) \geq 1 - \frac{\mathbb{E}\|G\|_{E \cap B_p^N}^p}{r^p} \geq 1 - O\left(\frac{dp^{p/2}}{r^p}\right).$$

Setting $K = E \cap rB_p^N$, then for some $r = O(\sqrt{p} \cdot d^{1/p})$, $\gamma_E(K) \geq \frac{1}{2}$, which concludes the proof by Banaszczyk's theorem [5]. $\square$

REMARK. The above estimate can actually be improved to give tail estimates as follows. Let $E$ be an $m$-dimensional subspace of $\mathbb{R}^n$. For $p > 2$ the function $x \mapsto \|x\|_p$ is Lipschitz with constant 1 on $\mathbb{R}^n$ and the Gaussian isoperimetric inequality shows that for every $\epsilon > 0$,

$$\gamma_E(E \cap (\mathbb{E}\|G\|_{E \cap B_p^n} + \epsilon)B_p^n) \geq 1 - e^{-\epsilon^2/2}.$$

Since $\mathbb{E}\|G\|_{E \cap B_p^n} \leq \mathbb{E}\|G\|_{B_p^k} \leq c\sqrt{p} \cdot m^{1/p}$ for some absolute constant $c$, then

$$\gamma_E(E \cap (c\sqrt{p} \cdot m^{1/p} + \epsilon)B_p^n) \geq 1 - e^{-\epsilon^2/2}.$$

3.5. *An application*: *covering numbers of convex hulls of points in $\ell_2$ by $B_p$ balls.* In this section, which is similar in spirit to the previous one, we use our results to give an infinite-dimensional extension of a classical inequality which bounds the minimal number of cubes $\varepsilon B_\infty^d$ required to cover a convex hull of a finite number of points in $\ell_2^d$ (this classical result depends on the maximum of $d$ and the number of points). Here, we are interested in finding upper bounds of the minimal number of cubes $\varepsilon B_\infty$ required to cover a convex hull of a finite number of points in $\ell_2$ depending only on $\varepsilon$ and the number of taken points. Since the structure of $\ell_\infty$ depends deeply on the chosen basis in $\ell_2$, a simple approximation argument is not enough to obtain our result.

The main result of this section, as described in the Introduction, is restated below:

PROPOSITION 26. *There exists an absolute constant $C > 0$ such that for every integer $m$, $\varepsilon > 0$ and $2 \leq p \leq \infty$, for all $x_1, \ldots, x_m$ in the unit ball of $\ell_2$,*

$$\log N(\mathrm{absconv}\{x_1, \ldots, x_m\}, \varepsilon B_p) \leq C \frac{\log m}{\varepsilon^{p/(p-1)}}.$$

PROOF. We first prove the proposition in the case when $p = \infty$. Since all $x_i$'s are in $B_2$ we can find an integer $d$ so that we can write $x_i = y_i + z_i$ with $y_i \in B_2^d$ and $\|z_i\|_\infty < \epsilon$ for all $i = 1, \ldots, m$. If the absolute convex hull of $y_1, \ldots, y_m$ can be covered by $N$ translates of $\epsilon B_\infty^d$, then the absolute convex hull of $x_1, \ldots, x_m$ can be covered by $N$ translates of $2\epsilon B_\infty$. So, it is enough to prove the result for the $y_i$'s.

Let $T : \ell_1^m \to \ell_2^d$ defined by $Te_i = y_i$ for all $i = 1, \ldots, m$, $E = \mathrm{span}\{y_1, \ldots, y_m\}$ and $G$ be a Gaussian vector in $E$. Since $\|x_i\|_2 \leq 1$, then by Sudakov's inequality [35],

$$\sup_{\varepsilon > 0} \varepsilon \sqrt{\log N(T(B_1^m), \varepsilon(B_2^d \cap E))} \leq \mathbb{E} \sup_{i=1,\ldots,m} |\langle G, y_i \rangle| \leq C\sqrt{\log m}.$$

Moreover, by the dual Sudakov inequality due to [26],

$$\sup_{\varepsilon > 0} \varepsilon \sqrt{\log N(B_2^d \cap E, \varepsilon B_\infty^d)} \leq \mathbb{E}\|G\|_{\ell_\infty^d \cap E},$$

and by Corollary 19, $\mathbb{E}\|G\|_{\ell_\infty^d \cap E} \leq \mathbb{E}\|G\|_{\ell_\infty^{\dim E}} \leq C\sqrt{\log m}$. Therefore,

$$\sup_{\varepsilon > 0} \varepsilon \sqrt{\log N(B_2^d \cap E, \varepsilon B_\infty^d)} \leq C\sqrt{\log m}.$$

Since the covering numbers are sub-additive,

$\log N(T(B_1^m), \varepsilon B_\infty^d)$

$$\leq \log N(T(B_1^m), \sqrt{\varepsilon}(B_2^d \cap E)) + \log N(\sqrt{\varepsilon}(B_2^d \cap E), \varepsilon B_\infty^d)$$
$$\leq C \cdot \frac{\log m}{\varepsilon}.$$

For a general $p \geq 2$, the proof follows by interpolation. Recall that for Banach spaces $X$, $Y$ and a compact operator $u\colon X \to Y$, the entropy numbers of $u$ are defined for every integer $k$ by

$$e_k(u\colon X \to Y) = \inf\{\varepsilon; N(u(B_X), \varepsilon B_Y) \leq 2^k\}.$$

Let $T$ be defined as before on $\ell_1^m$ by $Te_i = x_i$ for all $i = 1, \ldots, m$. It is well known (see Lemma 12.1.11 in [28]) that for every integer $k$,

$$e_{2k-1}(T\colon \ell_1^m \to \ell_p) \leq e_k(T\colon \ell_1^m \to \ell_2)^{2/p} e_k(T\colon \ell_1^m \to \ell_\infty)^{1-2/p}.$$

The above result for $p = \infty$, stated in terms of entropy numbers, is

$$e_k(T\colon \ell_1^m \to \ell_\infty) \leq C \cdot \frac{\log m}{k},$$

and Sudakov's inequality [35] is just

$$e_k(T\colon \ell_1^m \to \ell_2) \leq C \cdot \sqrt{\frac{\log m}{k}}.$$

Therefore,

$$e_{2k-1}(T\colon \ell_1^m \to \ell_p) \leq C \cdot \left(\frac{\log m}{k}\right)^{1-1/p},$$

as claimed. $\square$

**Acknowledgments.** We thank A. Giannopoulos and M. Talagrand for their stimulating questions.

F. BARTHE
LABORATOIRE DE STATISTIQUES
  ET PROBABILITÉS
CNRS UMR C5583
UNIVERSITÉ PAUL SABATIER
31062 TOULOUSE CEDEX 04
FRANCE
E-MAIL: barthe@math.ups-tlse.fr

O. GUÉDON
EQUIPE D'ANALYSE FONCTIONNELLE
  DE L'INSTITUT DE MATHÉMATIQUES
UNIVERSITÉ PARIS 6
PARIS 75005
FRANCE
E-MAIL: guedon@math.jussieu.fr

S. MENDELSON
CENTRE FOR MATHEMATICS
  AND ITS APPLICATIONS
AUSTRALIAN NATIONAL UNIVERSITY
CANBERRA, ACT 0200
AUSTRALIA
E-MAIL: shahar.mendelson@anu.edu.au

A. NAOR
THEORY GROUP
MICROSOFT RESEARCH
REDMOND, WASHINGTON 98052
USA
E-MAIL: anaor@microsoft.com